\documentclass{article}
\usepackage[%
journal=IFB,
lang=british,
]{ems-journal}

\usepackage{mathtools}
\usepackage{mathrsfs}

\theoremstyle{plain}
\newtheorem{theorem}{Theorem}[section]

\newtheorem{lemma}[theorem]{Lemma}

\theoremstyle{definition}
\newtheorem{definition}[theorem]{Definition}

\newtheorem{remark}[theorem]{Remark}

\numberwithin{equation}{section}

\begin{document}

\title{Correct mathematical models of joint filtration of two immiscible viscous liquids.}
\titlemark{Joint filtration of two liquids.}

\emsauthor{1}{
	\givenname{Anvarbek}
	\surname{Meirmanov}
	\mrid{}
	\zblid{}
	\orcid{0000-0002-8543-3897}}{A.~Meirmanov}

\Emsaffil{1}{
	\department{}
	\organisation{Institute of Physical and Technical Sciences, L.\,N. Gumilyov Eurasian National University}
	\rorid{}
	\address{Satbaev St., 2}
	\zip{010000}
	\city{Astana}
	\country{Kazakhstan}
	\affemail{anvarbey1947@gmail.com}}

\emsauthor{2}{
	\givenname{Akbota}
	\surname{Senkebayeva.}
	\mrid{}
	\zblid{}
	\orcid{0000-0001-6646-5415}}{A.~Senkebayeva}

\Emsaffil{2}{
	\department{}
	\organisation{School of Applied Mathematics, Kazakh-British Technical University}
	\rorid{}
	\address{Tole bi, 59}
	\zip{050000}
	\city{Almaty}
	\country{Kazakhstan}
	\affemail{akbota.senkebayeva@gmail.com}}    
\classification[35R35, 35R45]{35D30}


\begin{abstract}
Mathematical models of joint filtration of liquids are the main part of mathematical models of oil displacement by suspension. Since mining is a very important and urgent economic task, exact modeling of joint filtration of
two different fluids is also an urgent economic task. For example,
mathematical models of oil displacement by suspension are needed to create a hydrodynamic simulator of oil by suspension. All the existing simulators are based on the macroscopic Buckley--Leverett model, which does not distinguish between the free boundary separating liquids and the details of liquid interaction. All these fundamental processes occur at a microscopic level corresponding to the average size of pores, while all proposed macroscopic models operate on completely different orders of magnitude and do not distinguish between free boundaries or the characteristics of fluid interactions and are simply a set of axioms.
Exact modeling involves describing the process using the equations of classical Newtonian continuum mechanics at the microscopic level (average size of tens of micrometers). The only obstacle to using such models is that any numerical implementation for domains hundreds of meters in size would take years. A solution to this problem (the homogenization method) was proposed in the works of J. Keller and E. Sánchez-Palencia.
\end{abstract}

\maketitle

\section{Introduction.}
\label{sec1}
Mathematical models of joint filtration of various liquids are the main part of mathematical models of oil displacement by suspension or mathematical models of filtration of radioactive suspension from tailings storage facilities. Since mining and ecology are very important and urgent economic problems, exact modeling of joint filtration of two different fluids is also an important and urgent economic problem. It is well known that solid ground is a complex heterogeneous medium. Heterogeneity means that soil properties vary in space. Moreover, analyses of boreholes and cores show that physical properties of soil change in volumes of several tens of microns when both liquid and solid components of soil are present in these volumes. Such media are called heterogeneous porous media. The domain occupied by the fluid is called the pore space, and the domain complementing it is called the solid skeleton. In an ideal situation, an oil reservoir simulator is used to describe the displacement of oil by suspension. It should consist of a block containing a description of the structure of the solid skeleton, the physical characteristics of the solid skeleton, a block with mathematical models (at both the microscopic and macroscopic levels), and a block for visualising for these mathematical models with the corresponding initial and boundary conditions. 

Exact modeling involves describing this process using the equations of classical \textbf{\emph{Newtonian continuum mechanics}} in periodic media at the \textbf{\emph{Microscopic level}}. The validity of such a description has been confirmed over the centuries. The result of such modeling is \textbf{\emph{Exact microscopic mathematical models}} with a unique solution stable under small perturbations. This description is the best of all possible. The only obstacle to using such models at the microscopic level is that in domains of hundreds of meters in size, any numerical implementation of microscopic mathematical models (average size of tens of microns) with a spatial step of several microns will take years, which is completely unacceptable.

The main goal of our manuscript is both the practical significance of mathematical models describing the specified physical processes and the scientific significance of new mathematical problems and methods in modeling the physical processes under consideration. The problem studied here is a \emph{\textbf{Free-boundary problem}} describing the joint movement of two immiscible fluids in the pore space of a solid skeleton.

These problems can serve as a theoretical basis (\textbf{\emph{Prototype}}) for a hydrodynamic simulator of an oil reservoir. All these problems are highly relevant for many countries. The methods developed in our manuscript can be used to predict groundwater movement and groundwater contamination migration, and to assess aquifer stability (\emph{\textbf{Hydrogeology and Water Resources Management}}). The models obtained can be used to analyse the consequences of oil and industrial fluid spills, calculate contamination zones, and develop measures to localise them, which contributes to the development of environmental protection measures and improves environmental safety of industrial regions (\emph{\textbf{Environmental safety}}).

The prototype of any simulator is a complex consisting of mathematical models (both microscopic and macroscopic) comprising of a system of differential equations, supplemented by the corresponding boundary conditions of continuity of normal displacements of the medium and normal stresses of various components of the medium at the free (\emph{\textbf{Unknown}}) boundary $\Gamma^{\varepsilon}(t)$, separating different liquids, as well as boundary conditions at given boundaries and initial conditions.

The practical significance and the main objective of our investigation noted above consist in the creation of a prototype hydrodynamic simulator of an oil reservoir providing the most exact microscopic and macroscopic description of the physical processes of liquid filtration from an oil reservoir.

We will define a \textbf{\emph{Simulator of the physical process}} as a complex consisting of macroscopic models of this physical process for different structures of the solid skeleton, physical characteristics of the liquid and solid components of the studied medium, supplemented by programs for implementing these models.

It should be noted right away that there are several publications in the scientific literature on joint filtration of liquids \cite{1}, \cite{3}, \cite{4}, \cite{7}, \cite{13}, in which the behaviour of the components of the medium under consideration is described by a mathematical model consisting of differential equations and boundary and initial conditions. From the context of the articles, it is clear that the proposed models are the \textbf{\emph{Macroscopic phenomenological mathematical models}}, describing the process at the \textbf{\emph{Macroscopic level}}. Unlike \textbf{\emph{Microscopic Mathematical Models}}, with a characteristic size of approximately tens of microns, in \textbf{\emph{Macroscopic Mathematical Models}} the characteristic size is approximately decimeters or tens of decimeters. Because of this, these models do not distinguish between the microstructure of a continuous medium, since in such a model at each point the medium contains the solid skeleton, fluids in the pores of this skeleton and free  boundary, separating these liquids.

As for mathematical models of oil reservoirs, publications on this topic do not contain any conclusions about scientific and practical significance based on the analysis of the Buckley--Leverett model. The authors of the most recent publication \cite{13} rigorously analysed all possible structures of the solutions and proved that in the case of displacement of oil by a suspension there appears a \emph{\textbf{Mixed zone of oil and suspension}}, which is impossible for immiscible fluids.

The fundamental difference between the ideas presented in our article and previously published results lies in the initial formulation of the problems under consideration and the research methods, which will be described in detail below using the example of mathematical modelling of the filtration of immiscible liquids.

The Buckley--Leverett model is the only model traditionally considered as the basic mathematical model describing the process of oil displacement by a suspension for the hydrodynamic simulator \textbf{\emph{Eclipse}} from \textbf{\emph{Schlumberger}}, simulator \textbf{\emph{Tempest}} from \emph{\emph{Roxar}}, simulator \textbf{\emph{VIP}} from \textbf{\emph{Landmark}} and simulator \textbf{\emph{TimeZYX}} from \textbf{\emph{Standard Oil and Trust}}. Other well-known models are usually modifications or generalizations of the Buckley--Leverett model. Essentially, each of these models is based on a system of differential equations together with boundary and initial conditions, supplemented by equations of state that determine the physical laws of filtration. As a rule, fluid motion is described by \emph{\textbf{Darcy's law of filtration}} and \textbf{\emph{mass conservation laws}} for the concentration of each fluid, provided that the sum of these concentrations is equal to one. 

Obviously, it makes no sense to draw any conclusions from such publications, since these conclusions are only possible for \textbf{\emph{Exact macroscopic mathematical models}} obtained as an \textbf{\emph{Exact homogenization}} of exact microscopic mathematical models of a given physical process. All phenomenological models are a set of postulates by the authors of the model that has no relation to the physical process being described by J.~B.~Keller \cite{2} and E.~Sanchez-Palencia \cite{20}. They were the first to explain that the exact description at the macroscopic level is possible if and only if:

(a) the physical process under consideration is described at the microscopic level by equations of Newtonian classical continuum mechanics (\textbf{\emph{Exact mathematical models}});

(b) a set of small dimensionless parameters is selected;

(c) the macroscopic mathematical model is an exact asymptotic limit (\textbf{\emph{Homogenization}}) of exact mathematical models at the microscopic level, when the selected small parameters tend to zero.

In the proposed microscopic models, we assume that the solid skeleton has a periodic structure with a cell size of $\varepsilon$ and describe the physical process under consideration at the microscopic level using Newton's postulates of classical continuum mechanics. These are the Stokes equations for the liquid motion coupled with transport equations for the viscosities of the liquids in the pore space of an absolutely rigid solid skeleton. The last postulate means that the displacements and velocities of a rigid skeleton are identically equal to zero.

We use the notation adopted in \cite{9}, \cite{10}, \cite{12}.

\subsection{The problem statement.}

As we have already noted, the derivation of macroscopic mathematical models should be based on the most exact mathematical model of the physical process at the microscopic level, described by the laws of Newtonian classical continuum mechanics (see section~1.1, \S 1.2 in \cite{14}).

We will assume that the liquids in the pores are incompressible, and
the solid skeleton is an absolutely rigid body, that is,
\begin{equation}\label{eq1.1}
\boldsymbol{v}^{\varepsilon}_{s}(\boldsymbol{x},t)=\boldsymbol{w}^{\varepsilon}_{s}(\boldsymbol{x},t)=0,
\,\,\,(\boldsymbol{x},t)\in\Omega_{T}.
\end{equation}
The motion of liquids at $t>0$ in the domain $\Omega^{\varepsilon}_{sp}(t)$, occupied by the suspension, and the domain $\Omega^{\varepsilon}_{ol}(t)$, occupied by the oil, in dimensionless variables \cite{14} is described by the stationary Stokes equations for incompressible viscous liquids
\begin{equation}\label{eq1.2}
\nabla\cdot\mathbb{P}^{\varepsilon}_{j}=\nabla{p}^{0},\,\,\mathbb{P}^{\varepsilon}_{j}=
\mu_{j}\varepsilon^{2}\mathbb{D}(x,\boldsymbol{v}^{\varepsilon}_{j}),\,\,
\boldsymbol{v}^{\varepsilon}_{j}=\frac{\partial\boldsymbol{w}^{\varepsilon}_{j}}{\partial t},\,\,
\boldsymbol{x}\in\Omega^{\varepsilon}_{f,j}(t),\,\,j=ol, sp,
\end{equation}
transport equations
\begin{equation}\label{eq1.3}
\frac{d{\mu}_{j}}{dt}\equiv\frac{\partial{\mu}_{j}}{\partial{t}}+
(\nabla\mu_{j}\cdot\boldsymbol{v}^{\varepsilon}_{j})=0,\,\,
\boldsymbol{x}\in\Omega^{\varepsilon}_{f,j}(t),\,\,\mu_{j}(\boldsymbol{x},0)=\mu_{j}
\end{equation}
for viscosities of the liquids and linearized continuity equation
\begin{equation}\label{eq1.4}
\nabla\cdot\boldsymbol{v}^{\varepsilon}_{j}=\nabla\cdot\boldsymbol{w}^{\varepsilon}_{j}=0,\,\,
\boldsymbol{x}\in\Omega^{\varepsilon}_{f,j}(t),\,\,j=ol,\,sp
\end{equation}
for dynamic characteristics $\boldsymbol{w}^{\varepsilon}_{f}$ (\textbf{\emph{Liquid displacements}}),
$\boldsymbol{w}^{\varepsilon}_{sp}$ (\textbf{\emph{Displacements of the suspension}}),
$\boldsymbol{w}^{\varepsilon}_{ol}$ (\textbf{\emph{Displacements of the oil}}),
$\displaystyle\,\boldsymbol{v}^{\varepsilon}_{f}=
\frac{\partial\boldsymbol{w}^{\varepsilon}_{f}}{\partial t}$ (\textbf{\emph{Liquid velocity}}),
$\displaystyle\,\boldsymbol{v}^{\varepsilon}_{ol}=
\frac{\partial\boldsymbol{w}^{\varepsilon}_{ol}}{\partial t}$ (\textbf{\emph{Liquid velocity of the oil}}),
$\displaystyle\,\boldsymbol{v}^{\varepsilon}_{sp}=
\frac{\partial\boldsymbol{w}^{\varepsilon}_{sp}}{\partial t}$ (\textbf{\emph{Liquid velocity of the suspension}}),
${\mu}^{\varepsilon}_{ol}$ (\textbf{\emph{Oil viscosity}}) and ${\mu}^{\varepsilon}_{sp}$ (\textbf{\emph{Suspension viscosity}}).

The following boundary conditions 
\begin{equation}\label{eq1.5}
\boldsymbol{w}^{\varepsilon}_{sp}=\boldsymbol{w}^{\varepsilon}_{ol},\,\,
\boldsymbol{v}^{\varepsilon}_{sp}=\boldsymbol{v}^{\varepsilon}_{ol},\,\,
\mathbb{P}^{\varepsilon}_{sp}<\textbf{N}^{\varepsilon}>=
\mathbb{P}^{\varepsilon}_{ol}<\textbf{N}^{\varepsilon}>
\end{equation}
expressing the laws of conservation of momentum \cite{14} (see Appendix~A, section~A.7) are satisfied
at the \textbf{\emph{Free boundary}} $\Gamma^{\varepsilon}(t)$ between the liquid components in the pore space.

At the given boundaries with injection well ${S}^{1}$, production well ${S}^{2}$, impermeable boundary ${S}^{0}$ and at the internal boundary $S^{\varepsilon}$, separating pore space and solid skeleton, we assume that
\begin{equation}\label{eq1.6}
\begin{aligned}
&
\mathbb{P}^{\varepsilon}_{f}<\boldsymbol{n}>=0,\,\,
\boldsymbol{x}\in\,{S}^{\varepsilon}\cup{S}^{1}\cup{S}^{2},\,\,t>0,
\\
& \boldsymbol{w}^{\varepsilon}_{ol}(\boldsymbol{x},t)=
\boldsymbol{w}^{\varepsilon}_{sp}(\boldsymbol{x},t)=0,
\\
& \boldsymbol{v}^{\varepsilon}_{ol}(\boldsymbol{x},t)=
\boldsymbol{v}^{\varepsilon}_{sp}(\boldsymbol{x},t)=0,\,\,\boldsymbol{x}\in{S}^{0},\,\,t>0,
\\
& \boldsymbol{w}^{\varepsilon}_{ol}(\boldsymbol{x},t)=
\boldsymbol{w}^{\varepsilon}_{sp}(\boldsymbol{x},t)=0,\,\,\boldsymbol{x}\in{S}^{\varepsilon},\,\,t>0,
\end{aligned}
\end{equation}
where $\mathbb{P}^{\varepsilon}_{f}=\chi^{\varepsilon}_{sp}\mathbb{P}^{\varepsilon}_{sp}+
\chi^{\varepsilon}_{ol}\mathbb{P}^{\varepsilon}_{ol}$.

The problem is completed with the boundary condition
\begin{equation}\label{eq1.7}
{\mu}^{\varepsilon}_{sp}(x_{1},x_{2},-\tfrac{1}{2},t)={\mu}_{sp}
\end{equation}
at the boundary ${S}^{1}$ and the initial conditions
\begin{equation}\label{eq1.8}
{\mu}^{\varepsilon}_{ol}(\boldsymbol{x},0)={\mu}_{ol},\,\,
\boldsymbol{w}^{\varepsilon}_{j}(\boldsymbol{x},0)=0,\,j=ol,\,sp,\,\,\boldsymbol{x}\in\Omega,\,\,
\Gamma^{\varepsilon}(0)=\Gamma_{0}.
\end{equation}
Hereinafter, we will refer to the problem \eqref{eq1.1}--\eqref{eq1.8} as the problem $\mathbb{A}^{\varepsilon}$.

As an \textbf{\emph{Approximate problem}} $\mathbb{B}^{\varepsilon}$, we will call the problem $\mathbb{A}^{\varepsilon}$ with an additional term in the dynamic equations
\begin{equation}\label{eq1.9}
\nabla\cdot\mathbb{P}^{\varepsilon}_{j}-\varepsilon^{2}
\frac{\partial\boldsymbol{w}^{\varepsilon}_{j}}{\partial{t}}=\nabla{p}^{0},\,\,j=ol, sp.
\end{equation}

Let us take a quick look at the structure of the manuscript.

In section~\ref{sec2}, well-known facts and definitions are given, as well as results on the compactness of sequences in a periodic structure (Lemma~\ref{l2.2} and Lemma~\ref{l2.3}), results on the extension of functions defined in a pore space $\Omega^{\varepsilon}_{f,T}$ with a given structure onto the domain $\Omega^{\varepsilon}_{s,T}$ (Lemma~\ref{l2.13}) and orthogonal decomposition of the space $\mathbb{L}_{2}(\Omega)$ (Lemma~\ref{l2.9}). Next we formulate the integral identity for the problem $\mathbb{A}^{\varepsilon}$, equivalent to the corresponding system of differential equations supplemented with boundary and initial conditions, formally derive the final homogenised model $\mathbb{H}$ for joint filtration of two immiscible liquids and formulate the problem $\mathbb{B}^{\varepsilon}$, equivalent to the corresponding system of differential equations supplemented with boundary and initial conditions. 

At the end of this section, we shall derive a shock relation at the free boundary, which determines the normal velocity at that boundary in terms of dynamic characteristics.
 
In section~\ref{sec3}, we formulate the main results of the manuscript.

In section~\ref{sec4}, we prove the existence of a weak solution to the problem $\mathbb{B}^{\varepsilon}$.

Finally, in section~\ref{sec5}, we derive a homogenised model $\mathbb{H}$ for the joint filtration of two immiscible incompressible fluids.

The physical process we are considering has a rather long duration (the filtration rate of the liquid is several meters per year). Therefore, the most interesting are the theorems of the existence of solutions to the corresponding initial-boundary value problems globally in time. On the other hand, due to the strong nonlinearity of the \textbf{\emph{Free boundary problems}} \cite{15}, it is usually not possible to prove any result globally in time for mathematical models at the microscopic level. That is, the only possible results can be local in time theorems on the existence of a weak or classical solution to an initial boundary value problem for a system of differential equations describing a physical process under consideration \cite{15}.

As usual, almost every new problem has multiple choices for the problem setting. For example, in our case, we can consider the non-stationary Stokes equations, but then we somehow must find a priori estimates for the liquid velocities keeping in mind the difficulties with the free boundary separating liquid components.
In any case, we first need to find some approximations that will help simplify the problem. This is the path we have decided to take. But, in any case, the chosen approximation must allow us to somehow write down the approximate mathematical model in the form of a system of integral identities equivalent to the approximate system of differential equations, supplemented by the corresponding boundary conditions. 

These integral identities allow us to correctly determine weak solutions to the problem $\mathbb{B}^{\varepsilon}$ and find corresponding a'priori estimates.

We will obtain the required homogenised model $\mathbb{H}$ by taking the limit as $\varepsilon\rightarrow{0}$.

\section{Preliminaries.}
\label{sec2}

\subsection{Dimensionless parameters.}

The dimensionless parameter $\displaystyle\,\varepsilon=\frac{1}{n}$ is taken as a small parameter,
where $n$ is an integer.

The dimensionless parameters $\mu^{\varepsilon}_{j},\,\,j=ol,\,sp$ characterize the dynamics of
viscosities of the liquid in pores:
\begin{equation}\label{eq2.1}
\mu^{\varepsilon}_{j}=\mu_{j}\chi^{\varepsilon}_{j},\,\, 0<\mu_{*}\leqslant\mu_{j}\leqslant\mu_{*}^{-1}<\infty,\,\,j=ol,\,sp,
\end{equation}
where $\mu_{*}$ is a given positive constant, $\chi^{\varepsilon}_{j}(\boldsymbol{x},t)$ is the \textbf{\emph{Characteristic function}} of the domain $\Omega^{\varepsilon}_{j}(t),\,j=ol,\,sp$ and as subscripts $\textbf{\emph{ol}}$ and $\textbf{\emph{sp}}$ we denote physical characteristics of oil and suspension respectively.

The mathematical model \eqref{eq1.1}-\eqref{eq1.9} under condition \eqref{eq2.1} is called \textbf{\emph{Biot's model for joint filtration of immiscible liquids}}.

\subsection{Domains and boundaries.}

Let $\Omega\subset\mathbb{R}^{3}$ be a bounded domain with piecewise smooth boundary
$S=\partial\Omega={S}^{0}\cup{S}^{1}\cup{S}^{2}$, where as $\partial\Omega$
we denote the boundary of the domain $\Omega$,

$\Omega^{\varepsilon}_{s}$--\textbf{\emph{Solid skeleton}}, the domain occupied by the solid component,

$\Omega^{\varepsilon}_{f}$--\textbf{\emph{Pore space}}, the domain occupied by the liquids,

${S}^{\varepsilon}={\Omega}^{\varepsilon}_{s}\cap{\Omega}^{\varepsilon}_{f}$-boundary between domains $\Omega^{\varepsilon}_{s}$ and $\Omega^{\varepsilon}_{f}$,

The boundary $S^{0}\subset\partial\Omega$ is impermeable to liquid in the pore space,
the boundary $S^{1}\subset\partial\Omega$ simulates injection wells and the boundary
$S^{2}\subset\partial\Omega$ simulates production wells.

In what follows we will assume that $\displaystyle\,{\Omega}=\{\tfrac{1}{2}\,\leqslant\,x_{1},\,x_{2},\,x_{3}\,
\leqslant\,\tfrac{1}{2}\}$ is the unit cube, ${\Omega}_{T}={\Omega}\times(0,T)\subset\mathbb{R}^{4}$,

$\displaystyle\,S^{0}=\{\boldsymbol{x}:x_{3}=0,\,-\tfrac{1}{2}\,\leqslant\,x_{1},\,x_{2}\,
\leqslant\,\tfrac{1}{2}\}$,

$\displaystyle\,S^{1}=\{\boldsymbol{x}:x_{1}=
-\tfrac{1}{2},\,-\tfrac{1}{2}\,\leqslant\,x_{2},\,x_{3}\,\leqslant\,\tfrac{1}{2}\}$,

$S^{2}=\{\boldsymbol{x}:x_{1}=\tfrac{1}{2},\,-\tfrac{1}{2}\,
\leqslant\,x_{2},\,x_{3}\,\leqslant\,\tfrac{1}{2}\}$,

$\displaystyle\,S^{j}_{T}=S^{j}\times[0,T]$,

$\displaystyle\,\Gamma(0)=
\{\boldsymbol{x}\in\Omega:-\tfrac{1}{2}<x_{1},x_{2}<\tfrac{1}{2},\,x_{3}=--\tfrac{1}{2}\}$.

Next we set

$\Omega^{\varepsilon}_{sp}(t)\subset\Omega^{\varepsilon}_{f}$--domain in the pore space, occupied by the suspension,
$\Omega^{\varepsilon}_{ol}(t)\subset\Omega^{\varepsilon}_{f}$--domain in the pore space, occupied by the oil,
\begin{equation*}
\Omega^{\varepsilon}_{f,T}=\Omega^{\varepsilon}_{f}\times(0,T),\,\,
\Omega^{\varepsilon}_{s,T}=\Omega^{\varepsilon}_{s}\times(0,T),\,\,
\Gamma^{\varepsilon}_{T}=\bigcup_{t=0}^{t=T}\Gamma^{\varepsilon}(t),\,\,
\Omega^{\varepsilon}_{j,T}=\bigcup_{t=0}^{t=T}\Omega^{\varepsilon}_{j}(t),\,j=ol,\,sp
\end{equation*}
and
\begin{equation*}
\Omega_{f}^{\boldsymbol{k},\varepsilon}=\Omega^{\varepsilon}_{f}\cap{\Omega}^{\boldsymbol{k}},\,\,
\Omega_{s}^{\boldsymbol{k},\varepsilon}=
\Omega^{\varepsilon}_{s}\cap{\Omega}^{\boldsymbol{k},\varepsilon},\,\,
\Gamma^{\boldsymbol{k},\varepsilon}(t)=\Gamma^{\varepsilon}(t)\cap{\Omega}^{\boldsymbol{k}}
\end{equation*}
for all $\boldsymbol{k}=(k_{1},k_{2},k_{3})$, $k_{1},k_{2},k_{3}\in \mathbb{Z}$ (integer numbers).

For any continuous in $\Omega_{f,j}(t),\,j=sp,\,ol$ function $\boldsymbol{u}(\boldsymbol{x},t)$ its limits at the points $\boldsymbol{x}_{0}\in\Gamma^{\varepsilon}(t)$ are denoted as
\begin{equation*}
\begin{aligned}
& 
\boldsymbol{u}(\boldsymbol{x}_{0}+0,t)=
\lim_{\boldsymbol{x}\rightarrow \boldsymbol{x}_{0}}\boldsymbol{u}(\boldsymbol{x},t),
\,\,\boldsymbol{x}\in \Omega^{\varepsilon}_{ol}(t),\,\,\,
\boldsymbol{x}_{0}\in \Gamma^{\varepsilon}(t),
\\
& \boldsymbol{u}(\boldsymbol{x}_{0}-0,t)=\lim_{\boldsymbol{x}\rightarrow
\boldsymbol{x}_{0}}\boldsymbol{u}(\boldsymbol{x},t),\,\,
\boldsymbol{x}\in \Omega^{\varepsilon}_{sp}(t),\,\,\,\boldsymbol{x}_{0}\in \Gamma^{\varepsilon}(t).
\end{aligned}
\end{equation*}
As a small parameter we choose $\displaystyle\,\varepsilon=\frac{1}{n},\,\, n=1,2,3,\ldots$ so that the boundary condition \eqref{eq1.6} on the boundary ${S}^{1}\cup{S}^{2}$ makes sense.

\subsection{The structures of the solid skeleton and pore space.}

In what follows all functions of the type $\varphi(\boldsymbol{y};\boldsymbol{x})$, where $\boldsymbol{x}\in \Omega$ and $\boldsymbol{y}\in\textbf{Y}_{f}\subset\textbf{Y}=
\mathbb{R}^{3}$ are considered 1-periodic in variable $\boldsymbol{y}$:
\begin{equation*}
\varphi(\boldsymbol{y};x)=\varphi\big(\boldsymbol{\varsigma}(\boldsymbol{y}),x\big),\,\,\,
\boldsymbol{y}=[|\boldsymbol{y}|]+\varepsilon\,\boldsymbol{\varsigma}(\boldsymbol{y}),\,\,
[|\boldsymbol{y}|]=([|y_{1}|],\,[|y_{2}|],\,[|y_{3}|]).
\end{equation*}
The number $[|a|]$ is the integer part of the number $a$.

We restrict ourselves to the simplest structures of the solid skeleton $\textbf{Y}_{s}$ and pore space, where
\begin{equation*}
\begin{aligned}
& 
\textbf{Y}=\{\boldsymbol{y}\in\mathbb{R}^{3}:-\tfrac{1}{2}<y_{k}<\tfrac{1}{2},\,\,k=1,2,3\},
\\
& \textbf{Y}=\textbf{Y}_{f}\cup\gamma\cup\textbf{Y}_{s},
\,\,\textbf{Y}_{s}=\textbf{Y}\backslash\overline{\textbf{Y}}_{f},\,\,
\gamma=\overline{\textbf{Y}}_{f}\cap\overline{\textbf{Y}}_{s},
\\
& \textbf{Y}_{f}=\textbf{Y}^{(1)}_{f}\cup\textbf{Y}^{(2)}_{f}\cup\textbf{Y}^{(3)}_{f},
\{\boldsymbol{y}\in \textbf{Y}:-\tfrac{1}{2}<y_{1}<\tfrac{1}{2},\,|y|_{2}^{2}+|y_{3}|^{2}<r_{0}^{2}\},
\\
& \textbf{Y}^{(1)}_{f}=\{\boldsymbol{y}\in \textbf{Y}:
-\tfrac{1}{2}<y_{1}<\tfrac{1}{2},\,|y_{2}|^{2}+|y_{3}|^{2}<r_{0}^{2}\},
\\
& \textbf{Y}^{(2)}_{f}=\{\boldsymbol{y}\in \textbf{Y}:
-\tfrac{1}{2}<y_{2}<\tfrac{1}{2},\,|y_{1}|^{2}+|y_{3}|^{2}<r_{0}^{2}\},
\\
& \textbf{Y}^{(3)}_{f}=\{\boldsymbol{y}\in \textbf{Y}:-
\tfrac{1}{2}<y_{3}<\tfrac{1}{2},\,|y_{1}|^{2}+|y_{2}|^{2}<r_{0}^{2}\}
\end{aligned}
\end{equation*}
with a given function $r_{0}(\boldsymbol{x})$, $\displaystyle\,0<r_{0}(\boldsymbol{x})<\frac{1}{2}$.

As $\chi^{\varepsilon}_{f}(\boldsymbol{x})$ we denote the characteristic function of the pore space $\Omega^{\varepsilon}_{f}$: $\chi^{\varepsilon}_{f}(\boldsymbol{x})=1$ in the pore space $\Omega^{\varepsilon}_{f}$ and $\chi^{\varepsilon}_{f}(\boldsymbol{x})=0$ in the solid skeleton $\Omega^{\varepsilon}_{s}=\Omega\setminus\overline{\Omega}^{\,\varepsilon}_{f}$.

Let also $\chi^{\varepsilon}_{sp}(\boldsymbol{x},t)$ be the characteristic function of the liquid domain in the pore space occupied by the suspension and $\chi^{\varepsilon}_{ol}(\boldsymbol{x},t)$ be the characteristic function of the liquid domain in the pore space occupied by the oil.

Then by construction
\begin{equation}\label{eq.2.2}
\begin{aligned}
&
\chi^{\varepsilon}_{f}(\boldsymbol{x})=\chi\big(\boldsymbol{\varsigma}(\boldsymbol{y});
\boldsymbol{x}\big)=\chi\big(\boldsymbol{\varsigma}(\frac{\boldsymbol{x}}{\varepsilon});
\boldsymbol{x}),\,\,\chi^{\varepsilon}_{f}(\boldsymbol{x})=0,\,\,
\boldsymbol{x}\in\Omega^{\varepsilon}_{s},\,\,\chi^{\varepsilon}_{f}(\boldsymbol{x})=1,\,\,
\boldsymbol{x}\in\Omega^{\varepsilon}_{f},
\\
& \chi^{\varepsilon}_{sp}(\boldsymbol{x},t)=
\frac{1}{\mu_{sp}}\mu^{\varepsilon}_{sp}(\boldsymbol{x},t),\,\,
\chi^{\varepsilon}_{ol}(\boldsymbol{x},t)=
\frac{1}{\mu_{ol}}\mu^{\varepsilon}_{ol}(\boldsymbol{x},t),
\\
& \chi^{\varepsilon}_{sp}(\boldsymbol{x},t)=1,\,\,\boldsymbol{x}\in\Omega^{\varepsilon}_{f,sp}(t),\,\,
\chi^{\varepsilon}_{sp}(\boldsymbol{x},t)=0,\,\,\boldsymbol{x}\in\Omega^{\varepsilon}_{f,ol}(t),
\\
& \chi^{\varepsilon}_{ol}(\boldsymbol{x},t)=1,\,\,\boldsymbol{x}\in\Omega^{\varepsilon}_{f,ol}(t),\,\,
\chi^{\varepsilon}_{ol}(\boldsymbol{x},t)=0,\,\,\boldsymbol{x}\in\Omega^{\varepsilon}_{f,sp}(t),
\\
& \chi^{\varepsilon}_{ol}+\chi^{\varepsilon}_{sp}=\chi^{\varepsilon}_{f}.
\end{aligned}
\end{equation}

\subsection{Matrices, tensors and differential operators.}

We fix the standard Cartesian orthogonal basis
$\boldsymbol{e}^{1},\,\boldsymbol{e}^{2},\,\boldsymbol{e}^{3}$
in $\mathbb{R}^{3}$ and $\mathbb{A}$, $\mathbb{B}$ and $\mathbb{C}$ are \textbf{\emph{Tensors}} (linear transformations $\mathbb{R}^{3}\rightarrow\mathbb{R}^{3}$).
The action of the tensor $\mathbb{A}$ on the vector $\boldsymbol{b}$
is denoted as the vector $\boldsymbol{c}=\mathbb{A}<\boldsymbol{b}>$. As $(\boldsymbol{a}\cdot\boldsymbol{b})$ we denote the \textbf{\emph{Scalar product}} of vectors $\boldsymbol{a}$ and $\boldsymbol{b}$.
The product $\mathbb{C}=\mathbb{A}\cdot\mathbb{B}$ is a transformation
$\mathbb{R}^{3}\rightarrow\mathbb{R}^{3}$, where $\mathbb{A}, \mathbb{B}:\mathbb{R}^{3}\rightarrow\mathbb{R}^{3}$, 

$\mathbb{C}<\boldsymbol{x}>=\mathbb{A}<\mathbb{B}<\boldsymbol{x}>>$.

$\mathbb{I}$ is a unit tensor: $\mathbb{I}\cdot\mathbb{A}=\mathbb{A}\cdot\mathbb{I}=\mathbb{A}$ for any tensor $\mathbb{A}$ and $\mathbb{I}<\boldsymbol{a}>=\boldsymbol{a},\,\forall\boldsymbol{a}\in\mathbb{R}^{3}$.

For any vectors $\boldsymbol{a}$, $\boldsymbol{b}$, $\boldsymbol{c}$ as $\boldsymbol{a}\otimes\boldsymbol{b}$ we denote the \textbf{\emph{Diad}} (second-rank tensor), where 

$(\boldsymbol{a}\otimes\boldsymbol{b})<\boldsymbol{c}>=
\boldsymbol{a}(\boldsymbol{b}\cdot\boldsymbol{c})$.

As $\mathbb{J}^{ij}$ we denote the second-rank tensor $\mathbb{J}^{ij}=\frac{1}{2}
(\boldsymbol{e}^{i}\otimes\boldsymbol{e}^{j}+\boldsymbol{e}^{j}\otimes\boldsymbol{e}^{i})$.

Then $\displaystyle\,\mathbb{A}=\sum_{i,j=1}^{3}{a}_{ij}\boldsymbol{e}^{i}\otimes\boldsymbol{e}^{j}$.

The second order tensor $\mathbb{A}$ is symmetric, if
$(\mathbb{A}<\boldsymbol{e}^{j}>\cdot\boldsymbol{e}^{i})=
(\mathbb{A}<\boldsymbol{e}^{i}>\cdot\boldsymbol{e}^{j})$.

As is easy to see, any second-rank tensor is a linear mapping from $\mathbb{R}^{3}\rightarrow\mathbb{R}^{3}$.

By $(A)$, $(B)$ and $(C)$ we denote the corresponding to tensors $\mathbb{A}$, $\mathbb{B}$ and $\mathbb{C}$ matrices in the chosen Cartesian coordinate system:
\begin{equation*}
(A)=\left(
\begin{array}{ccc}
a_{11}&a_{12}&a_{13}\\
a_{21}&a_{22}&a_{23}\\
a_{31}&a_{32}&a_{33}
\end{array}
\right),\,\,(B)=
\left(
\begin{array}{ccc}
b_{11}&b_{12}&b_{13}\\
b_{21}&b_{22}&b_{23}\\
b_{31}&b_{32}&b_{33}
\end{array}
\right),\,\,(C)=
\left(
\begin{array}{ccc}
c_{11}&c_{12}&c_{13}\\
c_{21}&c_{22}&c_{23}\\
c_{31}&c_{32}&c_{33}
\end{array}
\right).
\end{equation*}
For matrices, the usual operations of sum $(A)+(B)$, multiplication by scalars $\alpha(A)$ and product $(B)\cdot(C)$ are defined.

Let $\widetilde{\boldsymbol{w}}(\boldsymbol{x},t)=
\big({w}_{1}(\boldsymbol{x},t),{w}_{2}(\boldsymbol{x},t),{w}_{3}(\boldsymbol{x},t)\big)$ and
$\mathbb{D}(x,\boldsymbol{w})=
\frac{1}{2}\big(\nabla_{x}\widetilde{\boldsymbol{w}}+(\nabla_{x}\widetilde{\boldsymbol{w}})^{*}\big)$.

Then the second-rank symmetric tensor
\begin{equation}\label{eq2.3}
\begin{aligned}
&
\mathbb{D}(x,\widetilde{\boldsymbol{w}})=\frac{1}{2}\sum_{i,j=1}^{3}({d}_{ij}(x,\widetilde{\boldsymbol{w}})
\boldsymbol{e}^{i}\otimes\boldsymbol{e}^{j}+
{d}_{ji}(x,\widetilde{\boldsymbol{w}})\boldsymbol{e}^{j}\otimes\boldsymbol{e}^{i}),
\\
& {d}_{ij}(x,\widetilde{\boldsymbol{w}})=
\frac{\partial\,w_{i}}{\partial\,x_{j}}(x,\widetilde{\boldsymbol{w}}),\,\,i,j=1,2,3,
\end{aligned}
\end{equation}
is called the \textbf{\emph{Symmetric gradient}} of the vector $\widetilde{\boldsymbol{w}}$.
By definition
\begin{equation}\label{eq2.4}
\begin{aligned}
&
\mathbb{D}(x,\widetilde{\boldsymbol{w}})<\boldsymbol{a}>\,\,\stackrel{df.}{=}\frac{1}{2}\sum_{i,j=1}^{3}
\big({d}_{ij}(x,\widetilde{\boldsymbol{w}})(\boldsymbol{e}^{i}\otimes
\boldsymbol{e}^{j})+{d}_{ji}(x,\widetilde{\boldsymbol{w}})
(\boldsymbol{e}^{j}\otimes\boldsymbol{e}^{i})\big)<\boldsymbol{a}>=
\\
& \frac{1}{2}\sum_{i=1}^{3}({d}_{ij}(x,\widetilde{\boldsymbol{w}}){a}_{j}\boldsymbol{e}^{i}+
{d}_{ji}(x,\widetilde{\boldsymbol{w}}){a}_{i}\boldsymbol{e}^{j}),
\\
& \mathbb{D}(x,\widetilde{\boldsymbol{w}}):\mathbb{D}(x,\widetilde{\boldsymbol{v}})=
\sum_{i,j=1}^{3}{d}_{ij}(x,\widetilde{\boldsymbol{w}}){d}_{ji}(x,\widetilde{\boldsymbol{v}}).
\end{aligned}
\end{equation}

\subsection{Strong convergence criteria in $\mathbb{L}_{2}(\Omega_{T}).$}

\begin{definition}\label{de1}
We say that the function ${\mu^{\varepsilon}_{j}}\in\mathbb{L}_{2}(\Omega_{T}),\,j=sp,\,ol$,
possesses a time derivative
$\displaystyle\,\frac{\partial\mu^{\varepsilon}_{j}}{\partial{t}}\in\mathbb{L}_{2}\big(0,T;
\mathbb{W}^{-1}_{2}(\Omega)\big)$, if
\begin{equation*}
\Big|\int_{0}^{T}\int_{\Omega}\mu^{\varepsilon}_{j}\frac{\partial \xi}{\partial t}dxdt\Big|\leqslant \,M_{\mu}\,\Big(\int_{0}^{T}\int_{\Omega}|\nabla\xi|^{2}dxdt\Big)^{\frac{1}{2}},\,j=sp,\,ol
\end{equation*}
with some positive constant $M_{\mu}$ independent of $\xi$ for the arbitrary functions $\xi\in\mathbb{W}^{1,1}_{2}(\Omega_{T})$.
\end{definition}

\begin{lemma}[\cite{12}]\label{l2.2}
Let the sequences $\{\mu^{\varepsilon}_{j}\}$ and $\{\boldsymbol{v}^{\varepsilon}_{j},\,j=sp,\,ol\}$
be uniformly bounded in the space $\mathbb{L}_{2}(\Omega_{T})$, and the sequence of derivatives
$\displaystyle\,\{\frac{\partial\mu^{\varepsilon}_{j}}{\partial{t}}\}$
be uniformly bounded in the space $\mathbb{L}_{2}\big(0,T;\mathbb{W}^{-1}_{2}(\Omega)\big)$.

Then there exists a subsequence of the sequence $\{\mu^{\varepsilon}_{j},\,j=sp,\,ol\}$, that converges strongly in $\mathbb{L}_{2}(\Omega_{T})$.
\end{lemma}

The generalization of this lemma for a periodic structure with characteristic function
$\displaystyle\,\chi^{\,\varepsilon}(\boldsymbol{x})=
\chi(\frac{\boldsymbol{x}}{\varepsilon};\boldsymbol{x})$
has been proved by A.~Meirmanov and O.~Galtsev \cite{16}:

\begin{lemma}[\cite{16}]\label{l2.3}
Let $\displaystyle\,\chi^{\varepsilon}(\boldsymbol{x})=
\chi(\frac{\boldsymbol{x}}{\varepsilon};\boldsymbol{x})$, where $\chi(\boldsymbol{y};\boldsymbol{x})$
be 1-periodic in $\boldsymbol{y}$ function, the sequences
$\{\mu^{\varepsilon}_{j}\}$ and $\{\boldsymbol{v}^{\varepsilon}_{j},\,j=sp,\,ol\}$ be uniformly bounded in $\mathbb{L}_{2}(\Omega_{T})$.

Then there exists some subsequence of $\{\mu^{\varepsilon}_{j}\},\,j=sp,\,ol$ that converges strongly in $\mathbb{L}_{2}(\Omega_{T})$.
\end{lemma}

\subsection{Two-scale convergent methods.}

\begin{definition}\label{de4}
\emph{The sequence $\{\widetilde{\boldsymbol{w}}^{\,\varepsilon}\}\subset \mathbb{L}_{2}(\Omega_{T})$,
is said to be two-scale convergent to the function
$\boldsymbol{W}(\boldsymbol{y};\boldsymbol{x},t)\in \mathbb{L}_{2}(\Omega_{T}\times Y)$,
which is 1-periodic in the variable $\boldsymbol{y}\in Y$ (notation
$\widetilde{\boldsymbol{w}}^{\,\varepsilon}\,\stackrel{\!\!2-sc.}{\rightarrow}
\,\boldsymbol{W}(\boldsymbol{y};\boldsymbol{x},t)$),
if for any smooth function $\sigma=\sigma(\boldsymbol{y};\boldsymbol{x},t)$,
1-periodic in the variable $\boldsymbol{y}$ the following equality holds}
\begin{equation}\label{eq2.5}
\begin{aligned}
&
\lim_{{\varepsilon}\to 0}\int\int_{\Omega_{T}}\widetilde{\boldsymbol{w}}^{\,\varepsilon}(\boldsymbol{x},t) \sigma(\frac{\boldsymbol{x}}{\varepsilon};\boldsymbol{x},t)dxdt=
\int\int_{\Omega_{T}}\big(\int_{\textbf{Y}}\boldsymbol{W}(\boldsymbol{y};\boldsymbol{x},t)
\sigma(\boldsymbol{y};\boldsymbol{x},t)dy\big)dxdt.
\end{aligned}
\end{equation}
\end{definition}
 us recall that the sequence $\{\widetilde{\boldsymbol{w}}^{\,\varepsilon}\}$ weakly converges to the function $\boldsymbol{w}$ if
\begin{equation*}
\lim_{{\varepsilon}\to 0}\int\int_{\Omega_{T}}
(\widetilde{\boldsymbol{w}}^{\,\varepsilon}\cdot\boldsymbol{\phi})dxdt=
\int\int_{\Omega_{T}}(\boldsymbol{w}\cdot\boldsymbol{\phi})dxdt
\end{equation*} for any smooth function $\boldsymbol{\phi}$.

Note that weak and two-scale convergence are connected by the relation:
\begin{equation}\label{eq2.6}
\begin{aligned}
&
\mbox{if}\,\,\widetilde{\boldsymbol{w}}^{\,\varepsilon}\,\stackrel{\!\! 2-sc.}{\rightarrow}\,
\boldsymbol{W}(\boldsymbol{y};\boldsymbol{x},t)\,\,\,\,\mbox{(converges two-scale)},
\\
& \mbox{then}\,\,\widetilde{\boldsymbol{w}}^{\,\varepsilon}(\boldsymbol{x},t)
\rightharpoonup \int_{\textbf{Y}}\boldsymbol{W}(\boldsymbol{y};\boldsymbol{x},t) d{y}=\boldsymbol{w}(\boldsymbol{x},t)
\,\,\mbox{(converges weakly)}.
\end{aligned}
\end{equation}
The existence and basic properties of two-scale convergent sequences are proved in the following theorem.

\begin{theorem}[\cite{17}]\label{t2.5}
Let the structure of the pore space be given by the function  $\displaystyle\,\chi(\frac{\boldsymbol{x}}{\varepsilon};\boldsymbol{x},t)$.

Then for any bounded in $\mathbb{W}^{1,0}_{2}(\Omega_{T})\cap
\mathbb{L}_{2}\big(0,T;\mathbb{W}^{-1}_{2}(\Omega)\big)$ sequence $\{\widetilde{\boldsymbol{w}}^{\varepsilon}\}$

\textbf{1.} Any bounded in $\mathbb{L}_{2}(\Omega_{T})$ sequence $\{\widetilde{\boldsymbol{w}}^{\,\varepsilon}\}$
contains some subsequence two-scale convergent to some function
$\boldsymbol{W}(\boldsymbol{y};\boldsymbol{x},t)$, $\boldsymbol{W}\in\mathbb{L}_{2}
(\Omega_{T}\times \textbf{Y})$, 1-periodic in the variable $\boldsymbol{y}$.

\textbf{2.} Let sequences $\{\widetilde{\boldsymbol{w}}^{\varepsilon}\}$ and $\{\varepsilon\mathbb{D}(x,\widetilde{\boldsymbol{w}}^{\,\varepsilon})\}$ be uniformly bounded in $\mathbb{L}_{2}(\Omega_{T})$.

Then there exists the function $\boldsymbol{W}=\boldsymbol{W}(\boldsymbol{y};\boldsymbol{x},t)$, 1-periodic
in $\boldsymbol{y}$, and the subsequence $\{\widetilde{\boldsymbol{w}}^{\,\varepsilon}\}$ (for simplicity we will keep the same notations) such that
$\boldsymbol{W},\,\mathbb{D}(y,\boldsymbol{W})\in \mathbb{L}_{2}(\Omega_{T}\times\textbf{Y})$,
and sequences $\{\widetilde{\boldsymbol{w}}^{\,\varepsilon}\}$ and $\{\varepsilon\mathbb{D}(x,\widetilde{\boldsymbol{w}}^{\,\varepsilon})\}$ (for simplicity we keep the same indices for subsequences) two-scale converge in
$\mathbb{L}_{2}(\Omega_{T})$ to $\boldsymbol{W}$ and $\mathbb{D}(y,\boldsymbol{W})$ correspondingly.

\textbf{3.} Let sequences $\{\widetilde{\boldsymbol{w}}^{\,\varepsilon}\}$ and
$\{D(x,\widetilde{\boldsymbol{w}}^{\,\varepsilon})\}$ be bounded in $\mathbb{L}_{2}(\Omega_{T})$.

Then there are functions $\boldsymbol{w}(\boldsymbol{x},t)\in \mathbb{W}^{1,0}_{2}(\Omega_{T})$ and
$\boldsymbol{W}(\boldsymbol{y};\boldsymbol{x},t),\,\boldsymbol{W}\in \mathbb{L}_{2}(\Omega_{T}\times \textbf{Y})\cap\mathbb{W}^{1,0}_{2}(\textbf{Y})$,
subsequence of $\{\mathbb{D}(x,\widetilde{\boldsymbol{w}}^{\,\varepsilon})\}$ such that the function
$\boldsymbol{W}$ is 1-periodic in $\boldsymbol{y}$, $\mathbb{D}(x,\boldsymbol{w}) \in \mathbb{L}_{2}(\Omega_{T})$,
$D(y,\boldsymbol{W}) \in \mathbb{L}_{2}(\Omega_{T}\times\textbf{Y})$, and the sequence
$\{\mathbb{D}(x,\widetilde{\boldsymbol{w}}^{\varepsilon})\}$ two-scale converges to the function
$\mathbb{D}(x,\boldsymbol{w})+D(y,\boldsymbol{W})$.
\end{theorem}

\subsection{Poincar\'{e} inequality.}

\begin{lemma}[\cite{19}]\label{l2.6}
Let $\Omega\subset\mathbb{R}^{3}$ be a bounded domain with a piecewise smooth
Lipschitz boundary. Then for any function
$\widetilde{\boldsymbol{w}}^{\varepsilon}\in\stackrel{\!\!\circ}{\mathbb{W}}^{1}_{2}(\Omega)$,
\begin{equation*}
\|\widetilde{\boldsymbol{w}}^{\varepsilon}_{f}\|_{2,\Omega}
\leqslant\varepsilon{M}_{\Omega}
\|\mathbb{D}(x,\widetilde{\boldsymbol{w}}^{\varepsilon})\|_{2,\Omega}.
\end{equation*}
If $\displaystyle\,\Omega\subset\bigcup_{|\boldsymbol{k}|=1}^{n^{3}}\Omega^{\boldsymbol{k},\varepsilon}$
and $\widetilde{\boldsymbol{w}}^{\varepsilon}\in\stackrel{\!\!\circ}{\mathbb{W}}^{1}_{2}
(\Omega^{\,\boldsymbol{k},\varepsilon})$,
$\boldsymbol{k}=(k_{1},k_{2},k_{3})\in\mathbb{Z}^{3}$, then
\begin{equation*}
\int_{\Omega^{\boldsymbol{k},\varepsilon}}
|\widetilde{\boldsymbol{w}}^{\varepsilon}(\boldsymbol{x})|^{2}dx
\leqslant\varepsilon^{2}\,M_{\Omega}
\int_{\Omega^{\boldsymbol{k},\varepsilon}}
|\mathbb{D}\big(x,\widetilde{\boldsymbol{w}}^{\varepsilon}(\boldsymbol{x})\big)|^{2}dx
\end{equation*}
and
\begin{equation}\label{eq2.7}
\int_{\Omega}|\widetilde{\boldsymbol{w}}^{\varepsilon}(\boldsymbol{x})|^{2}dx
\leqslant\varepsilon^{2}\,M_{\Omega}\int_{\Omega}
|\mathbb{D}\big(x,\widetilde{\boldsymbol{w}}^{\varepsilon}(\boldsymbol{x})\big)|^{2}dx.
\end{equation}
The constant $M_{\Omega}$ is bounded for bounded domain $\Omega$.
\end{lemma}

\begin{remark}\label{2.7}
A similar result
\begin{equation*}
\int_{\textbf{Y}_{f}}(|\widetilde{\boldsymbol{w}}^{\varepsilon}-
\widetilde{\boldsymbol{w}}^{\varepsilon}_{\Omega}|^{2}dx\,\leqslant\,\varepsilon^{2}\,M_{\Omega}
\int_{\Omega}|\mathbb{D}
\big(x,\widetilde{\boldsymbol{w}}^{\varepsilon}(\boldsymbol{x})\big)|^{2}dx,
\end{equation*}
where $\displaystyle\,\widetilde{\boldsymbol{w}}^{\varepsilon}_{\Omega}=
\frac{1}{|\Omega|}\int_{\Omega}\widetilde{\boldsymbol{u}}^{\varepsilon}(\boldsymbol{x})dx$ and $|\Omega|$
is a volume of the domain $\Omega$, is called the \textbf{\emph{Poincar\'{e}--Wirtinger inequality}} for any $\widetilde{\boldsymbol{u}}^{\varepsilon}\in\mathbb{W}^{1}_{2}(\Omega)$.

For the proof, see \cite{8}.
\end{remark}

\subsection{The simplest embedding lemma.}

\begin{lemma}\label{l2.8}
Let $\Omega\subset\mathbb{R}^{3}$ with a piecewise smooth boundary from $\mathbb{C}^{1}$.

Then for any function $\widetilde{\boldsymbol{w}}_{f}\in\mathbb{W}^{1}_{2}(\Omega)$ identically equal to zero on some part of the boundary $\partial\Omega$ with strictly positive surface measure the following estimate
\begin{equation}\label{eq2.8}
\|\widetilde{\boldsymbol{w}}_{f}\|_{2,\Omega}\leqslant{M}
\|\mathbb{D}(x,\widetilde{\boldsymbol{w}}_{f})\|_{2,\Omega}
\end{equation}
holds true.

The constant ${M}$ is bounded for bounded $\Omega$.
\end{lemma}
For details, see \cite{5}, \cite{21}.

\subsection{Orthogonal decomposition of the space $\mathbb{L}_{2}(\Omega).$}

Let $\Omega\subset\mathbb{R}^{3}$ and $\stackrel{\,\,\!\!\circ}{\mathbb{J}}(\overline{\Omega})$ be the set of all infinitely smooth solenoidal vector functions $\boldsymbol{w}^{J}$ in $\Omega$, and $\stackrel{\!\!\circ}{\mathbb{G}}(\Omega)$ be the set of the gradients $\boldsymbol{w}^{G}=\nabla\,{w}^{G}$ of all infinitely smooth real functions ${w}^{G}$ in $\Omega$.

As $\stackrel{\!\!\circ}{\mathbb{J}^{(1)}}(\overline{\Omega})$ we denote the closure of the set $\stackrel{\!\!\circ}{\mathbb{J}}(\overline{\Omega})$ in $\mathbb{L}_{2}(\Omega)$ and as the set $\stackrel{\!\!\circ}{\mathbb{G}^{(1)}}(\overline{\Omega})$ we denote the closure of the set $\stackrel{\!\!\circ}{\mathbb{G}}(\overline{\Omega})$ in $\mathbb{L}_{2}(\Omega)$.

\begin{lemma}\label{l2.9}
The space $\mathbb{L}_{2}(\Omega)$ is the direct sum of the subspaces $\stackrel{\!\!\circ}{\mathbb{J}^{(1)}}({\overline\Omega})$ and $\stackrel{\!\!\circ}{\mathbb{G}^{(1)}}(\overline{\Omega})$ in $\mathbb{L}_{2}(\Omega)$.

In particular, for every solenoidal function $\boldsymbol{w}$ there exists a scalar function 
${p}\in\mathbb{W}^{1}_{2}(\Omega)$ such that $\boldsymbol{w}=\nabla{p}$ and 
\begin{equation*}
\|\nabla{p}\|_{2,\Omega}\leqslant{M}\|\boldsymbol{w}\|^{(1)}_{2,\Omega},
\end{equation*}
where constant $M$ does not depend on $\boldsymbol{w}$.
\end{lemma}
For the proof, see \cite{11}, chapter~1, $\S\,2$ and Lemma

\subsection{H\"{o}lder's inequality.} 

\begin{lemma}\label{l2.10}
For any $\widetilde{\boldsymbol{w}},\,\widetilde{\boldsymbol{v}}\in\mathbb{L}_{2}(\Omega)$ the H\"{o}lder's inequality  
\begin{equation}\label{eq2.9}
\|\widetilde{\boldsymbol{w}}\,\widetilde{\boldsymbol{v}}\|_{1,\Omega}\leqslant
\|\widetilde{\boldsymbol{w}}\|_{2,\Omega}\|\widetilde{\boldsymbol{v}}\|_{2,\Omega}
\end{equation}
holds.
\end{lemma}
For details, see \cite{5}, \cite{21}.

\subsection{Some functional spaces.}

By $\mathbb{L}_{2}(\Omega)$ we denote the functional space of all measurable functions
$u(\boldsymbol{x})$ with a finite norm
\begin{equation*}
\|u\|_{2,\Omega}=\big(\int_{\Omega}u^{2}(\boldsymbol{x})dx\big)^{\frac{1}{2}}.
\end{equation*}
The functional space $\mathbb{W}^{1}_{2}(\Omega)$ is the closure of all infinitely smooth in $\Omega$ functions in the norm
\begin{equation*}
\|u\|^{(1)}_{2,\Omega}=\|u\|_{2,\Omega}+\|\nabla\,u\|_{2,\Omega}.
\end{equation*}
As a space $\stackrel{\!\!\circ}{\mathbb{W}}^{1,0}_{2}(\Omega_{T})$ we define the space of all functions from
$\mathbb{W}^{1,0}_{2}(\Omega_{T})$ vanishing at the boundary $\partial{\Omega}$.

The functional space $\mathbb{C}^{k}(\overline{\Omega})$ for integer $k\geqslant 0$ consists of all
functions $u(\boldsymbol{x})$ with a finite norm
\begin{equation*}
|u|^{(k)}_{\Omega}=\max_{\boldsymbol{x}\in \Omega}\sum_{|m|=0}^{k}|D^{m}u(\boldsymbol{x})|,\,\,
D^{m}u=\frac{\partial^{|m|}u}{\partial x_{1}^{m_{1}}\cdots\partial x_{n}^{m_{n}}},
\end{equation*}
where $m=(m_{1},\ldots,m_{n}),\,m_{i}\geqslant 0,\,i=1,\ldots,n$, is a multi-index,
$|m|=m_{1}+\cdots+m_{n}$.

$\stackrel{\,\!\!\circ}{\mathbb{C}}^{\,k}(\overline{\Omega})$ is a subspace of
$\mathbb{C}^{k}(\overline{\Omega})$ of all functions $u(\boldsymbol{x})$ vanishing at the boundary $S=\partial\Omega$.

As a space $\mathbb{C}^{\infty}(\overline{\Omega})$ we denote the space of all infinitely smooth in $\Omega$ functions $u(\boldsymbol{x})$.

As a space $\mathbb{H}^{\alpha}(\overline{\Omega})$ we define the space of all functions $u(\boldsymbol{x})$ with a finite norm
\begin{equation*}
|u|^{(\alpha)}_{\Omega}=|u|^{(0)}_{\Omega_{T}}+\max_{\boldsymbol{x}\in\Omega,\,|h|<\infty}
\frac{|u(\boldsymbol{x}+h)-u(\boldsymbol{x})|}{|h|^{\alpha}}<\infty,\,\,
|u|^{(0)}_{\Omega_{T}}=\max_{(\boldsymbol{x},t)\in\Omega_{T}}|u(\boldsymbol{x},t)|.
\end{equation*}
As a space $\mathbb{H}^{\alpha,\frac{\alpha}{2}}(\overline{\Omega}_{t_{0}})$
we define the space of all functions $u(\boldsymbol{x},t)$ with a finite norm
\begin{equation*}
|u|^{(\alpha,\frac{\alpha}{2})}_{\Omega_{T}}=|u|^{(0)}_{\Omega_{T}}+
\max_{(\boldsymbol{x},t)\in\Omega_{T},\,|h|<\infty}
\frac{|u(\boldsymbol{x}+h,t+\frac{h}{2})-u(\boldsymbol{x},t)|}{|h|^{\alpha}}<\infty.
\end{equation*}
Finally, as a space $\mathbb{H}^{k+\alpha,\frac{k+\alpha}{2}}(\overline{\Omega}_{t_{0}})$ we define the space of all functions $u(\boldsymbol{x},t)$ with a finite norm
\begin{equation*}
|u|^{(k+\alpha,\frac{k+\alpha}{2})}_{\Omega_{T}}=|u|^{(0)}_{\Omega_{T}}+
\max_{(\boldsymbol{x},t)\in\Omega_{T},\,|h|<\infty}
\frac{|D^{k}u(\boldsymbol{x}+h,t+\frac{h}{2})-D^{k}u(\boldsymbol{x},t)|}{|h|^{\alpha}}<\infty.
\end{equation*}
Let $\mathbb{W}$ be some functional space with elements $\boldsymbol{w}(\boldsymbol{x},t)$.
Then as $\mathbb{L}^{\infty}\big(0,T;\mathbb{W}(\Omega)\big)$ we denote all functions $\boldsymbol{w}$, bounded in $\mathbb{W}(\Omega)$.

In particular, the space $\mathbb{L}^{\infty}\big(0,T;\mathbb{H}^{2+\alpha}(\overline{\Omega})\big)$ consists of all bounded in $\mathbb{H}^{2+\alpha}(\overline{\Omega})$ functions.

\subsection{Equivalent formulation of the problem $\mathbb{A}^{\varepsilon}$ as an integral identity.}

\begin{definition}\label{d2.11}
Let
$\mathbb{P}^{\varepsilon}_{ol}=\chi^{\varepsilon}_{ol}\varepsilon^{2}\mu_{ol}
\mathbb{D}(x,\widetilde{\boldsymbol{v}}_{ol}^{\varepsilon})$,
\,\,$\mathbb{P}^{\varepsilon}_{sp}=\varepsilon^{2}\mu_{sp}
\mathbb{D}(x,\boldsymbol{v}_{sp}^{\varepsilon})$.

We say that functions
$\boldsymbol{v}^{\varepsilon}_{j}
\in\mathbb{W}^{1,0}_{2}\big(\Omega_{j,T}\big)\,j=sp,\,ol$ are a weak solution to the problem $\mathbb{A}^{\varepsilon}$
if hold boundary and initial conditions \eqref{eq1.7}, \eqref{eq1.8} and the integral identities
\begin{equation}\label{eq2.10}
\begin{aligned}
&
-\int_{0}^{t_{0}}\int_{\Omega}\mathbb{P}^{\varepsilon}_{f}:\mathbb{D}(x,\boldsymbol{\varphi})dxdt=
\int_{0}^{t_{0}}(\nabla{p}^{0}\cdot\boldsymbol{\varphi})dxdt,
\\
& \int_{0}^{t_{0}}\int_{\Omega}\big((\chi^{\varepsilon}_{sp}
\boldsymbol{v}^{\varepsilon}_{sp}+\chi^{\varepsilon}_{ol}
\boldsymbol{v}^{\varepsilon}_{ol})\cdot\boldsymbol{\varphi}\big)dxdt=0
\end{aligned}
\end{equation}
for any arbitrary smooth solenoidal functions $\boldsymbol{\varphi}$ satisfying the boundary conditions
\begin{equation}\label{eq2.11}
\begin{aligned}
&
\boldsymbol{\varphi}(\boldsymbol{x}_{0}+0)=\boldsymbol{\varphi}(\boldsymbol{x}_{0}-0),\,\,
\boldsymbol{x}_{0}\in\Gamma^{\varepsilon}(t),
\\
& \boldsymbol{\varphi}(\boldsymbol{x}_{0}+0)=
\lim_{\boldsymbol{x}\rightarrow\boldsymbol{x}_{0}}\boldsymbol{\varphi}(\boldsymbol{x}),
\,\,\boldsymbol{x}\in\Omega^{\varepsilon}_{ol}(t_{0}),\,\,
\boldsymbol{x}_{0}\in\Gamma^{\varepsilon}(t),
\\
& \boldsymbol{\varphi}(\boldsymbol{x}_{0}-0)=
\lim_{\boldsymbol{x}\rightarrow\boldsymbol{x}_{0}}
\boldsymbol{\varphi}(\boldsymbol{x}),\,\,\boldsymbol{x}\in \Omega^{\varepsilon}_{sp}(t_{0}),\,\,
\boldsymbol{x}_{0}\in \Gamma^{\varepsilon}(t),
\\
& \boldsymbol{\varphi}(\boldsymbol{x},0)=\boldsymbol{\varphi}(\boldsymbol{x},T)=0,\,\,
\boldsymbol{\varphi}(\boldsymbol{x},t)=0,\,\,\,(\boldsymbol{x},t)
\in\big({S}^{\varepsilon}\cup{S}^{1}\cup{S}^{2}\big)\times(0,T).
\end{aligned}
\end{equation}
In \eqref{eq2.10}, \eqref{eq2.11} $0<t_{0}\leqslant{T}$.
\end{definition}

\subsection{Formal homogenisation of the problem $\mathbb{A}^{\varepsilon}.$}

As \begin{equation}\label{eq2.12}
{\pi}_{j}(\boldsymbol{x},t)=\int_{0}^{t}{p}_{j}(\boldsymbol{x},\tau)d\tau,\,\,j=sp,\,ol
\end{equation}
we denote an antiderivative of the liquid pressures ${p_{j}},\,\,j=sp,\,ol$.  

\begin{lemma}\label{l2.12}
Under the conditions of Theorem~\ref{t3.1} the formal homogenisation $\mathbb{H}$ of the problem $\mathbb{A}^{\varepsilon}$ consists of Darcy's law of filtration
\begin{equation}\label{eq2.13}
\boldsymbol{w}_{j}=-\frac{1}{\mu_{j}}(B)<\nabla(\pi_{j}-{p}^{0}t)>,\,\,
\nabla\cdot\boldsymbol{w}_{j}=0
\end{equation}
for the displacements $\boldsymbol{w}_{j}$ and the antiderivative ${\pi}_{j}$, \emph{j=ol,sp}
in the domain $\Omega_{T}$.

Differential equations are completed with the boundary and initial conditions
\begin{equation}\label{eq2.14}
{\pi}_{j}(\boldsymbol{x},t)={p}^{0}(\boldsymbol{x})t,\,j=sp,ol,\,\,
\boldsymbol{x}\in{S}^{1}\cup{S}^{2},\,\,\,0<t<T,
\end{equation}
\begin{equation}\label{eq2.15}
\boldsymbol{w}_{j}\cdot\boldsymbol{n}=0,\,j=sp,ol,\,\,\boldsymbol{x}\in{S}^{0},\,\,\,0<t<T.
\end{equation}
In \eqref{eq2.14} $\boldsymbol{n}$ is a normal vector to the boundary ${S}^{1}\cup{S}^{2}$, a symmetric strictly positive definite constant matrix $(B)$ is defined by the formula \eqref{eq5.16}.
\end{lemma}

\subsection{Extension Lemma.}

Extension results are very important in homogenization (\cite{6}, \cite{22}).

For example, some sequence has different properties in different domains and only the properties of the sequence in the first domain permit choosing a convergent subsequence. Therefore, we must preserve the best properties of the sequence and apply the extension from the first domain onto the second one.

\begin{lemma}\label{l2.13}
Let $\{\boldsymbol{w}^{\varepsilon}_{j}\}$ be a bounded sequence in
$\mathbb{W}^{1,0}_{2}(\Omega^{\varepsilon}_{j,T}),\,j=ol,\,sp$.

Then for all $\varepsilon>0$ there exists an extension operator

$\mathbb{E}_{j}:\mathbb{W}^{1,0}_{2}\big(\Omega_{j,T}\big)\rightarrow\mathbb{W}^{1,0}_{2}(\Omega_{T})$,
$\mathbb{E}_{j}({\boldsymbol{w}}_{j}^{\varepsilon})=
\widetilde{\boldsymbol{w}}^{\varepsilon}_{j}\in\mathbb{W}^{1,0}_{2}(\Omega_{T}),\,j=ol,\,sp$, such that
\begin{equation}\label{eq2.16}
\begin{aligned}
&
\varepsilon^{2}\int_{0}^{t_{0}}\int_{\Omega}\big(\chi^{\varepsilon}_{sp}\mu_{sp}
\mathbb{D}(x,\widetilde{\boldsymbol{v}}^{\varepsilon}_{sp})+\chi^{\varepsilon}_{ol}\mu_{ol}
\mathbb{D}(x,\widetilde{\boldsymbol{v}}^{\varepsilon}_{ol})\big):
\mathbb{D}(x,\boldsymbol{\varphi})dxdt=
\\
& -\varepsilon^{2}\int_{0}^{t_{0}}\int_{\Omega}\big(\chi^{\varepsilon}_{sp}\mu_{sp}
\mathbb{D}(x,\widetilde{\boldsymbol{w}}^{\varepsilon}_{sp})+\chi^{\varepsilon}_{ol}\mu_{ol}
\mathbb{D}(x,\widetilde{\boldsymbol{w}}^{\varepsilon}_{ol})\big):
\mathbb{D}(x,\frac{\partial\boldsymbol{\varphi}}{\partial{t}})dxdt=
\\
& -\int_{0}^{t_{0}}\int_{\Omega}\chi^{\varepsilon}_{f}\Big(\big(\nabla({p}^{0}t)+\varepsilon^{2}
\chi^{\varepsilon}_{sp}\widetilde{\boldsymbol{w}}^{\varepsilon}_{sp}+\chi^{\varepsilon}_{ol}
\varepsilon^{2}\widetilde{\boldsymbol{w}}^{\varepsilon}_{ol}\big)
\cdot\frac{\partial\boldsymbol{\varphi}}{\partial{t}}\Big)dxdt,
\\
& \int_{0}^{t_{0}}\int_{\Omega}\big((\chi^{\varepsilon}_{sp}\widetilde{\boldsymbol{w}}^{\varepsilon}_{sp}+
\chi^{\varepsilon}_{ol}\widetilde{\boldsymbol{w}}^{\varepsilon}_{ol})\cdot
\frac{\partial\boldsymbol{\varphi}}{\partial{t}}\big)dxdt=0.
\end{aligned}
\end{equation}
\end{lemma}
To prove this lemma we just put
\begin{equation}\label{eq2.17}
\begin{aligned}
&
\widetilde{\boldsymbol{w}}^{\varepsilon}_{ol}=\boldsymbol{w}^{\varepsilon}_{ol},\,\,
\mathbb{D}(x,\widetilde{\boldsymbol{w}}_{ol}^{\varepsilon})=
\mathbb{D}(x,\boldsymbol{w}_{ol}^{\varepsilon}),\,\,\boldsymbol{x}\in\Omega^{\varepsilon}_{ol}(t),
\\
& \widetilde{\boldsymbol{w}}^{\varepsilon}_{ol}=\boldsymbol{w}^{\varepsilon}_{sp},\,\,
\mathbb{D}(x,\widetilde{\boldsymbol{w}}_{ol}^{\varepsilon})=
\mathbb{D}(x,\boldsymbol{w}_{sp}^{\varepsilon}),\,\,\boldsymbol{x}\in\Omega^{\varepsilon}_{sp}(t),
\\
& \widetilde{\boldsymbol{w}}^{\varepsilon}_{sp}=\boldsymbol{w}^{\varepsilon}_{ol},\,\,
\mathbb{D}(x,\widetilde{\boldsymbol{w}}_{sp}^{\varepsilon})=
\mathbb{D}(x,\boldsymbol{w}_{ol}^{\varepsilon}),\,\,\boldsymbol{x}\in\Omega^{\varepsilon}_{ol}(t),
\\
& \widetilde{\boldsymbol{w}}^{\varepsilon}_{sp}=\boldsymbol{w}^{\varepsilon}_{sp},\,\,
\mathbb{D}(x,\widetilde{\boldsymbol{w}}_{sp}^{\varepsilon})=
\mathbb{D}(x,\boldsymbol{w}_{sp}^{\varepsilon}),\,\,\boldsymbol{x}\in\Omega^{\varepsilon}_{sp}(t),
\end{aligned}
\end{equation}
\begin{equation*}
\begin{aligned}
&
\widetilde{\boldsymbol{v}}^{\varepsilon}_{ol}=\boldsymbol{v}^{\varepsilon}_{ol},\,\,
\mathbb{D}(x,\widetilde{\boldsymbol{v}}_{ol}^{\varepsilon})=
\mathbb{D}(x,\boldsymbol{v}_{ol}^{\varepsilon}),\,\,\boldsymbol{x}\in\Omega^{\varepsilon}_{ol}(t),
\\
& \widetilde{\boldsymbol{v}}^{\varepsilon}_{ol}=\boldsymbol{v}^{\varepsilon}_{sp},\,\,
\mathbb{D}(x,\widetilde{\boldsymbol{v}}_{ol}^{\varepsilon})=
\mathbb{D}(x,\boldsymbol{v}_{sp}^{\varepsilon}),\,\,\boldsymbol{x}\in\Omega^{\varepsilon}_{sp}(t),
\\
& \widetilde{\boldsymbol{v}}^{\varepsilon}_{sp}=\boldsymbol{v}^{\varepsilon}_{ol},\,\,
\mathbb{D}(x,\widetilde{\boldsymbol{v}}_{sp}^{\varepsilon})=
\mathbb{D}(x,\boldsymbol{v}_{ol}^{\varepsilon}),\,\,\boldsymbol{x}\in\Omega^{\varepsilon}_{ol}(t),
\\
& \widetilde{\boldsymbol{v}}^{\varepsilon}_{sp}=\boldsymbol{v}^{\varepsilon}_{sp},\,\,
\mathbb{D}(x,\widetilde{\boldsymbol{v}}_{sp}^{\varepsilon})=
\mathbb{D}(x,\boldsymbol{v}_{sp}^{\varepsilon}),\,\,\boldsymbol{x}\in\Omega^{\varepsilon}_{sp}(t).
\end{aligned}
\end{equation*}

\subsection{Equivalent formulation of the problem $\mathbb{B}^{\varepsilon}$ as an integral identity.}

\begin{definition}\label{d2.14}
Let $\mathbb{P}^{\varepsilon}=\chi^{\varepsilon}_{sp}\mathbb{P}^{\varepsilon}_{sp}+
\chi^{\varepsilon}_{ol}\mathbb{P}^{\varepsilon}_{ol}$,
$\mathbb{P}^{\varepsilon}_{ol}=\chi^{\varepsilon}_{ol}\varepsilon^{2}\mu_{ol}
\mathbb{D}(x,\boldsymbol{v}^{\varepsilon}_{ol})$,
$\mathbb{P}^{\varepsilon}_{sp}=\chi^{\varepsilon}_{sp}\varepsilon^{2}\mu_{sp}
\mathbb{D}(x,\boldsymbol{v}^{\varepsilon}_{sp})$. 

Let also $\boldsymbol{\varphi}$ be arbitrary
smooth solenoidal functions that vanish at the boundary $\big(S^{1}\cup S^{2}\big)\times(0,T)$,
at $t=0,\,T$ and satisfying conditions \eqref{eq2.11} and $\eta$ be arbitrary smooth functions that
vanish at the boundary $\partial\Omega_{T}$.

We say that functions $\boldsymbol{v}^{\varepsilon}_{j}\in\mathbb{W}^{1,0}_{2}(\Omega_{j,T})$ 
and $\mu_{j}\in\mathbb{L}_{\infty}(\Omega_{T})\,\,j=sp,\,ol$ are a weak solution to the problem $\mathbb{B}^{\varepsilon}$, if the following integral identities
\begin{equation}\label{eq2.18}
\begin{aligned}
&
0=\varepsilon^{2}\int_{0}^{t_{0}}\int_{\Omega}\big(\chi^{\varepsilon}_{sp}\mu_{sp}
\mathbb{D}(x,\widetilde{\boldsymbol{v}}^{\varepsilon}_{sp})+
\chi^{\varepsilon}_{ol}\mu_{ol}\mathbb{D}(x,\boldsymbol{v}^{\varepsilon}_{ol})+
({p}^{\varepsilon}-{p}^{\varepsilon})\mathbb{I})\big):\mathbb{D}(x,\boldsymbol{\varphi})dxdt+
\\
& \int_{0}^{t_{0}}\int_{\Omega}\big(\nabla{p}^{0}+\varepsilon^{2}(\chi^{\varepsilon}_{sp}
\widetilde{\boldsymbol{v}}^{\varepsilon}_{sp}+\chi^{\varepsilon}_{ol}
\widetilde{\boldsymbol{v}}^{\varepsilon}_{ol})\cdot\boldsymbol{\varphi}\big)dxdt=0,
\\
& \int_{0}^{t_{0}}\int_{\Omega}\big(\frac{\partial\eta}{\partial{t}}
(\chi^{\varepsilon}_{sp}\mu_{sp}+\chi^{\varepsilon}_{ol}\mu_{ol})+
(\mu_{sp}\widetilde{\boldsymbol{v}}^{\varepsilon}_{sp}+
\mu_{ol}\widetilde{\boldsymbol{v}}^{\varepsilon}_{ol})\cdot\nabla\eta\big)dxdt=0,
\end{aligned}
\end{equation}
\begin{equation}\label{eq2.19}
\int_{0}^{t_{0}}\int_{\Omega}\big((\chi^{\varepsilon}_{sp}
\boldsymbol{v}^{\varepsilon}_{sp}+\chi^{\varepsilon}_{ol}
\boldsymbol{v}^{\varepsilon}_{ol})\cdot\boldsymbol{\varphi}\big)dxdt=0
\end{equation}
together with the boundary and initial conditions \eqref{eq1.7}, \eqref{eq1.8}.
\end{definition}

\begin{remark}\label{2.15}
In the first identity \eqref{eq2.18}, we have taken into account Lemma~\ref{l2.9}, which implies the equality $\big(\nabla(p^{\varepsilon}-p^{0})\cdot\boldsymbol{\varphi}\big)=0$.
\end{remark}

\begin{remark}\label{2.16}
For the homogenization procedure we will also use identities
\begin{equation}\label{eq2.20}
\begin{aligned}
& 
\varepsilon^{2}\int_{0}^{t_{0}}\int_{\Omega}\big(\chi^{\varepsilon}_{sp}\mu_{sp}
\mathbb{D}(x,\boldsymbol{w}^{\varepsilon}_{sp})+\chi^{\varepsilon}_{ol}\mu_{ol}
\mathbb{D}(x,\boldsymbol{w}^{\varepsilon}_{ol})\big):
\mathbb{D}(x,\frac{\partial\boldsymbol{\varphi}}{\partial{t}})dxdt=
\\
& -\int_{0}^{t_{0}}\int_{\Omega}\Big(\big(\nabla({p}^{0}t)+\varepsilon^{2}
\chi^{\varepsilon}_{sp}\boldsymbol{w}^{\varepsilon}_{sp}+\chi^{\varepsilon}_{ol}
\varepsilon^{2}\boldsymbol{w}^{\varepsilon}_{ol}\big)
\cdot\frac{\partial\boldsymbol{\varphi}}{\partial{t}}\Big)dxdt,
\end{aligned}
\end{equation}
\begin{equation}\label{eq2.21}
\int_{0}^{t_{0}}\int_{\Omega}\big((\chi^{\varepsilon}_{sp}
\boldsymbol{w}^{\varepsilon}_{sp}+\chi^{\varepsilon}_{ol}
\boldsymbol{w}^{\varepsilon}_{ol})\cdot\frac{\partial\boldsymbol{\varphi}}{\partial{t}}\big)dxdt=0,
\end{equation}
which are equivalent to the first identity in \eqref{eq2.18} and identity \eqref{eq2.19} under the same conditions \eqref{eq2.11} for the functions $\boldsymbol{\varphi}$.
\end{remark}

\subsection{Shock relations.}

As a basic mathematical model we consider the model $\mathbb{M}_{2}$ in the form of abstract equation (A.6.3) in \cite{14}, Appendix~A, section~A.6
\begin{equation}\label{eq2.22}
\frac{\partial}{\partial{t}}(\rho{u})+\nabla\cdot(\rho{u}\boldsymbol{v})=0,
\end{equation}
where we put $\textbf{X}=0$ and $Y$=0.

Then the shock relation at the free boundary $\Pi$ for this case in takes the form
\begin{equation}\label{eq2.23}
[\rho{u}(v_{N}-V_{N})]=0,\,\,\,\boldsymbol{x}\in\Pi,\,\,0<t<T.
\end{equation}
Here $\boldsymbol{v}$ is a velocity of the liquid, $\boldsymbol{N}$ is unit the normal vector to the free boundary $\Pi$, $v_{N}$ is the normal velocity of the liquid at the free boundary, $V_{N}$ is the normal velocity of the free boundary $\Pi$. 

Coming back to our notations we obtain: 
\begin{equation}\label{eq2.24}
\rho{u}=\mu^{\varepsilon}_{j}=\chi^{\varepsilon}_{j}\mu_{j},\,\,
\boldsymbol{v}=\boldsymbol{v}^{\varepsilon}_{j},\,\,
v_{N}={V}^{\varepsilon}_{N},\,\,j=ol,\,sp.
\end{equation}
\begin{lemma}\label{l2.17}
(Integration by parts \cite{14}, Appendix~A, section~A.6, formula~(A.6.13).)

Let $\boldsymbol{v}^{\varepsilon}_{j}$ and ${\mu}^{\varepsilon}_{j},\,j=sp,ol$ be a weak solution to the problem $\mathbb{B}^{\varepsilon}$, $\mathbb{C}^{1}$ boundary $\Gamma^{\varepsilon}(t)$ divides $\Omega^{\varepsilon}_{f}$ into two subdomains $\Omega^{\varepsilon}_{sp,T}$ and $\Omega^{\varepsilon}_{ol,T}$, $[{\mu}]={\mu}_{sp}-{\mu}_{ol}$ be the jump of ${\mu}$ at the boundary $\Gamma^{\varepsilon}(t)$. Then  
\begin{equation}\label{eq2.25}
{V}^{\varepsilon}_{N}={v}^{\varepsilon}_{N},\,\,\boldsymbol{x}\in\Gamma^{\varepsilon}(t).
\end{equation}
\end{lemma}

\begin{remark}\label{4}
The boundary condition \eqref{eq2.25} corresponds to the first boundary condition in the system (A.6.14), section~A.6, \cite{14}) with the equal densities.
\end{remark}

\begin{remark}\label{2.18}
For more details see Appendix~A, section~A.6 in \cite{14}.
\end{remark}

\section{Main results.}
\label{sec3}

\begin{theorem}\label{t3.1}
Let $\displaystyle\,0<r_{0}<\frac{1}{2}$, \, ${p}^{0}\in\mathbb{H}^{2+\alpha}(\overline{\Omega})$ and
$\displaystyle\,\frac{\partial{p^{0}}}{\partial{x}_{3}}(\boldsymbol{x})\,\leqslant\,0$
for $\boldsymbol{x}\in\Omega$.

Then under conditions \eqref{eq2.1} the problem $\mathbb{B}^{\varepsilon}$ has a unique weak solution
\begin{equation}\label{eq3.1}
\boldsymbol{w}^{\varepsilon}_{j}\in\mathbb{W}^{1,0}_{2}(\Omega_{j,T}),\,\boldsymbol{v}^{\varepsilon}_{j}
\in\mathbb{W}^{1,0}_{2}(\Omega_{j,T}),\,j=sp,\,ol.
\end{equation}
\end{theorem}

\begin{theorem}\label{t3.2}
Under the conditions of Theorem~\ref{t3.1} the problem $\mathbb{H}$ has a unique solution
\begin{equation*}
\displaystyle\,\boldsymbol{w}_{j}(.,t)\in\mathbb{L}^{\infty}(0,T;\mathbb{H}^{2+\alpha,\frac{2+\alpha}{2}}
\big(\overline{\Omega}_{T})\big),\,\,{p}_{j}(.,t)\in\mathbb{L}^{\infty}(0,T;
\mathbb{H}^{1+\alpha,\frac{1+\alpha}{2}}\big(\overline{\Omega}_{T})\big),\,\,j=sp,\,ol
\end{equation*}
and
\begin{equation}\label{eq3.2}
\max_{0<t<T}\big(|\boldsymbol{w}_{j}(.,t)|^{(2+\alpha)}_{\Omega_{T}}+
|{p}_{j}(.,t)|^{(1+\alpha)}_{\Omega_{T}}\big)\leqslant\,M\sqrt{T}M_{\Omega},\,\,j=sp,\,ol,
\end{equation}
where $M_{\Omega}$ is defined in Lemma~\ref{l2.6}.
\end{theorem}

\section{Proof of Theorem~\ref{t3.1}: existence of a weak solution to the problem $\mathbb{B}^{\varepsilon}.$}
\label{sec4}

\begin{lemma}\label{l4.1}
Under the conditions of Theorem~\ref{t3.1} the problem $\mathbb{B}^{\varepsilon}$ has a unique
weak solution such that
\begin{equation}\label{eq4.1}
\begin{aligned}
& 
\|\boldsymbol{v}^{\varepsilon}_{ol}\|_{2,\Omega_{T}}+\|\boldsymbol{v}^{\varepsilon}_{sp}\|_{2,\Omega_{T}}+
\varepsilon\|\boldsymbol{v}^{\varepsilon}_{f}\|_{2,\Omega_{T}}+
\\
& \varepsilon\|\mathbb{D}(x,\boldsymbol{v}^{\varepsilon}_{ol})\|_{2,\Omega_{T}}+
\varepsilon\|\mathbb{D}(x,\boldsymbol{v}^{\varepsilon}_{sp})\|_{2,\Omega_{T}}
\leqslant\,{M}\sqrt{T}{M}_{\Omega}\|\nabla{p}^{0}\|_{2,\Omega},
\\
& \max_{0<t<T}\Big(\|\boldsymbol{w}^{\varepsilon}_{ol}(.,t)\|_{2,\Omega}+
\|\boldsymbol{w}^{\varepsilon}_{sp}(.,t)\|_{2,\Omega}+
\varepsilon\|\mathbb{D}\big(x,\boldsymbol{w}^{\varepsilon}_{ol}(.,t)\big)\|_{2,\Omega}+
\varepsilon\|\mathbb{D}\big(x,\boldsymbol{w}^{\varepsilon}_{sp}(.,t)\big)\|_{2,\Omega}\Big) 
\\
& \leqslant\,{M}\sqrt{T}{M}_{\Omega}\|\nabla{p}^{0}\|_{2,\Omega},\,\,j=ol,\,sp,
\end{aligned}
\end{equation}
where $M$ does not depend on $\varepsilon$, $\Omega$ and ${p}^{0}$.
\end{lemma}

\begin{proof}
To estimate liquid velocities and displacements we put in \eqref{eq2.18} $\boldsymbol{\varphi}=\boldsymbol{v}^{\varepsilon}_{f}$:
\begin{equation}\label{eq4.2}
\begin{aligned}
& 
\int_{0}^{t_{0}}\int_{\Omega^{\varepsilon}_{f}}\big((\chi^{\varepsilon}_{sp}
\boldsymbol{v}^{\varepsilon}_{sp}+\chi^{\varepsilon}_{ol}\boldsymbol{v}^{\varepsilon}_{ol})
\cdot\nabla{p}^{0}\big)dxdt=
\\
& \varepsilon^{2}\int_{0}^{t_{0}}\int_{\Omega^{\varepsilon}_{f}}\big(\chi^{\varepsilon}_{sp}
\mu_{sp}\mathbb{D}(x,\boldsymbol{v}^{\varepsilon}_{sp})+\chi^{\varepsilon}_{ol}
\mu_{ol}\mathbb{D}(x,\boldsymbol{v}^{\varepsilon}_{ol})\big):
\mathbb{D}(x,\boldsymbol{v}^{\varepsilon}_{f})dxdt=
\\
& \varepsilon^{2}\int_{0}^{t_{0}}\int_{\Omega^{\varepsilon}_{f}}\big(\chi^{\varepsilon}_{sp}\mu_{sp}
|\mathbb{D}(x,\boldsymbol{v}^{\varepsilon}_{sp})|^{2}+\mu_{ol}\chi^{\varepsilon}_{ol}
|\mathbb{D}(x,\boldsymbol{v}^{\varepsilon}_{ol})|^{2}\big)dxdt=
\\
& \varepsilon^{2}\int_{0}^{t_{0}}\int_{\Omega^{\varepsilon}_{f}}\big(\mu_{sp}
|\mathbb{D}(x,\boldsymbol{v}^{\varepsilon}_{sp})|^{2}+\mu_{ol}
|\mathbb{D}(x,\boldsymbol{v}^{\varepsilon}_{ol})|^{2}\big)dxdt.
\end{aligned}
\end{equation}
Next applying Poincar\'{e} and H\"{o}lder's inequalities \eqref{eq2.7} and \eqref{eq2.9}, as well as the trivial inequality $\displaystyle\,{a}\,{b}\leqslant\frac{\delta}{2}{a}^{2}+\frac{{b}^{2}}{2\delta}$, we will obtain
\begin{equation}\label{eq4.3}
\begin{aligned}
& 
\mu_{*}\varepsilon^{2}\int_{\Omega^{\varepsilon}_{f}}
\big(|\mathbb{D}(x,\boldsymbol{v}^{\varepsilon}_{sp})|^{2}+
|\mathbb{D}(x,\boldsymbol{v}^{\varepsilon}_{ol})|^{2}\big)dx\leqslant
\\
& \varepsilon^{2}\int_{\Omega^{\varepsilon}_{f}}
\big(\mu_{sp}|\mathbb{D}(x,\boldsymbol{v}^{\varepsilon}_{sp})|^{2}+\mu_{ol}
|\mathbb{D}(x,\boldsymbol{v}^{\varepsilon}_{ol})|^{2}\big)dx,
\\
& \Big|\int_{\Omega^{\varepsilon}_{f}}\Big(\chi^{\varepsilon}_{sp}
(\boldsymbol{v}^{\varepsilon}_{sp}\cdot\nabla{p}^{0})+
\chi^{\varepsilon}_{ol}(\boldsymbol{v}^{\varepsilon}_{ol}\cdot\nabla{p}^{0})\Big)dxdt\Big|\leqslant
\\
& \frac{\delta}{2}\int_{0}^{t_{0}}\int_{\Omega^{\varepsilon}_{f}}
(|\boldsymbol{v}^{\varepsilon}_{ol}|^{2}+|\boldsymbol{v}^{\varepsilon}_{sp}|^{2})dx+
\frac{T}{2\delta}\int_{\Omega^{\varepsilon}_{f}}|\nabla{p}^{0}|^{2}dx\leqslant
\\
& \frac{\delta}{2}\varepsilon^{2}{M}_{\Omega}\int_{\Omega^{\varepsilon}_{f}}
\big(|\mathbb{D}(x,\boldsymbol{v}^{\varepsilon}_{sp})|^{2}+
|\mathbb{D}(x,\boldsymbol{v}^{\varepsilon}_{ol})|^{2}\big)dx+
\frac{T}{2\delta}\int_{\Omega^{\varepsilon}_{f}}|\nabla{p}^{0}|^{2}dx.
\end{aligned}
\end{equation}
The desired estimates for the velocities follow from the last inequality for
$\displaystyle\,\delta=\frac{2\mu_{*}}{{M}_{\Omega}}$.

To estimate the displacements, we use representations
\begin{equation*}
\boldsymbol{w}^{\varepsilon}_{j}(\boldsymbol{x},t)=
\int_{0}^{t}\boldsymbol{v}^{\varepsilon}_{f}(\boldsymbol{x},\tau)d\tau,\,\,
\mathbb{D}\big(x,\boldsymbol{v}^{\varepsilon}_{f}(\boldsymbol{x},t)\big)=
\int_{0}^{t}\mathbb{D}\big(x,\boldsymbol{w}^{\varepsilon}_{f}(\boldsymbol{x},\tau)\big)d\tau,\,\,j=ol,\,sp,
\end{equation*}
which prove the desired estimates \eqref{eq4.1}.
\end{proof}

\section{Proof of Theorem~\ref{t3.2}: homogenization of the problem $\mathbb{B}^{\varepsilon}.$}
\label{sec5}

The homogenization procedure itself is well explained in many publications (see \cite{6}, \cite{22} and references therein)

\begin{remark}\label{5.1}
We have reconstructed the pressures ${p}_{ol}$ and ${p}_{sp}$ and their antiderivatives ${\pi}_{ol}$ and ${\pi}_{sp}$ by decomposing the space $\mathbb{L}_{2}(\Omega)$ into a direct sum of the subspace of all
solenoidal vector functions and the subspace of all gradients of scalar functions (see Lemma~\ref{l2.9}).
\end{remark}
\begin{lemma}\label{l5.2}
Under the conditions of Theorem~\ref{t3.1} there exist unique functions $\boldsymbol{w}_{ol}$, $\boldsymbol{w}_{sp}$, $p_{ol}$, $p_{sp}$, $\displaystyle\,{\pi}_{ol}=\int_{0}^{t}p_{ol}(\boldsymbol{x},\tau)d\tau$ and $\displaystyle\,{\pi}_{sp}\int_{0}^{t}p_{sp}(\boldsymbol{x},\tau)d\tau$ such that
$p_{ol}, p_{sp},\boldsymbol{w}_{ol}, \boldsymbol{w}_{sp}\in\mathbb{W}^{1,0}_{2}(\Omega_{T})$,
$\boldsymbol{W}_{ol}, \boldsymbol{W}_{sp}\in\mathbb{L}_{2}\big(0,T;\mathbb{W}^{1}_{2}(\textbf{Y}_{f})\big)$.

1) The sequence $\{\widetilde{\boldsymbol{w}}^{\,\varepsilon}_{j}\}$ converges weakly to the function $\boldsymbol{w}_{j}$ and two-scale to the function
$\boldsymbol{W}_{j}(\boldsymbol{y};\boldsymbol{x},t),\,j=ol,\,sp$.

2) The sequences $\{\varepsilon\nabla\cdot\widetilde{\boldsymbol{w}}^{\,\varepsilon}_{j}\}$ converge two-scale to the functions $\nabla_{y}\cdot\boldsymbol{W}_{j},\,j=ol,\,sp$.

3) The sequences $\{\varepsilon\mathbb{D}(x,\widetilde{\boldsymbol{w}}^{\,\varepsilon}_{j})\}$ converge two-scale to the functions $\mathbb{D}(y,\boldsymbol{W}_{j}),\,j=ol,\,sp$.

4) The following a priori estimates hold true
\begin{multline}\label{eq5.1}
\|{\boldsymbol{w}}_{j}\|_{2,\Omega_{T}}+\|\boldsymbol{W}_{j}\|_{2,\textbf{Y}_{f}\times\Omega_{T}}+
\|\nabla_{y}\boldsymbol{W}_{j}\|_{2,{\textbf{Y}}_{f}\times\Omega_{T}}+
\|\mathbb{D}\big(y,\boldsymbol{W}_{j})\|_{2,{Y}_{f}\times\Omega_{T}}\leqslant
\\
\,{M}\sqrt{T}\|\nabla{p}^{0}\|_{2,\Omega},\,\,j=ol,\,sp,
\end{multline}
\begin{equation}\label{eq5.2}
\|{p}_{j}\|_{2,\Omega_{T}}+\|\nabla{p}_{j}\|_{2,\Omega_{T}}+\|{\pi}_{j}\|_{2,\Omega_{T}}+
\Big\|\nabla\frac{\partial{\pi}_{j}}{\partial{t}}\Big\|_{2,\Omega_{T}}\leqslant\,
{M}\sqrt{T}\|\nabla{p}^{0}\|_{2,\Omega},\,\,j=ol,\,sp,
\end{equation}
where the constant $M$ is bounded for bounded domains $\Omega$.
\end{lemma}
The proof is straightforward and based on the estimates \eqref{eq4.1}, Theorem~\ref{t2.5}, section~\ref{sec2}, Lemma~\ref{l5.2} and Remark~\ref{4}. For details see Lemma 5.2 in \cite{15}.

We only note that
\begin{equation*}
\varepsilon^{2}\int_{0}^{t_{0}}\int_{\Omega}\chi^{\varepsilon}
\Big|\frac{\partial\widetilde{\boldsymbol{w}}^{\,\varepsilon}_{f}}{\partial{t}}\Big|^{2}dxdt
\leqslant\,{M}\sqrt{T}\|\nabla{p}^{0}\|_{2,\Omega}
\end{equation*}
and
\begin{equation*}
\lim_{\varepsilon\rightarrow{0}}\varepsilon\int_{0}^{t_{0}}\int_{\Omega}\chi^{\varepsilon}
\Big|\frac{\partial\widetilde{\boldsymbol{w}}^{\,\varepsilon}_{f}}{\partial{t}}\Big|dxdt=0.
\end{equation*}

\begin{lemma}\label{l5.3}
Under the conditions of Theorem~\ref{t3.1} the limiting procedure in the integral identities \eqref{eq2.19} results in the following dynamic problem $\mathbb{H}$ for displacements and pressures of the liquid components consisting of Darcy's law of filtration
\begin{equation}\label{eq5.3}
\boldsymbol{w}_{j}=\int_{\textbf{Y}_{f}}\boldsymbol{W}_{j}dy=
-\frac{1}{\mu_{j}}(B)<\nabla_{x}({\pi}_{j}-{p}^{0}t)>,\,\,
\nabla_{x}\cdot\boldsymbol{w}_{j}=0,\,\,j=ol,\,sp
\end{equation}
for the liquid displacements $\boldsymbol{w}_{j}$ and the antiderivative ${\pi}_{j}$ of the pressure ${p}_{j},\,\,j=ol,\,sp$ in the domain $\Omega_{T}$, completed with the boundary conditions
\begin{equation}\label{eq5.4}
{\pi}_{j}(\boldsymbol{x},t)={p}^{0}t,\,\,j=ol,\,sp,\,\,\boldsymbol{x}\in{S}^{1}\cup{S}^{2},\,\,
\boldsymbol{w}_{j}\cdot\boldsymbol{n}=0,\,\,j=ol,\,sp,\,\,\boldsymbol{x}\in{S}^{0},\,\,0<t<T,
\end{equation}
where $\boldsymbol{n}$ is a normal vector to the boundary ${S}^{0}$.

The symmetric, strictly positive-definite constant matrix $(B)$ is given by the formula \eqref{eq5.16}.
\end{lemma}

\begin{proof}
First, we derive the continuity equations for the functions $\boldsymbol{w}_{j}$ and $\boldsymbol{W}_{j},\,j=ol,\,sp.$

To do this, let us consider the second integral identity in \eqref{eq2.22}
\begin{equation*}
0=\lim_{\varepsilon\rightarrow{0}}\int_{0}^{{t}_{0}}
\int_{\Omega}{\eta}(\nabla\cdot\boldsymbol{w}^{\varepsilon}_{j})dxdt=
-\int_{0}^{{t}_{0}}\int_{\Omega}(\nabla{\eta}\cdot\boldsymbol{w}_{j})dxdt,\,\,\,j=ol,\,sp.
\end{equation*}
This identity gives us continuity equations and boundary conditions for the displacements of oil and suspension
\begin{equation}\label{eq5.5}
\nabla\cdot\boldsymbol{w}_{j}=0,\,\,\boldsymbol{x}\in\Omega,\,\,(\boldsymbol{w}_{j}\cdot\boldsymbol{n})=0, \,\,\boldsymbol{x}\in{S}^{0},\,\,j=ol,\,sp.
\end{equation}
To derive the continuity equation for the unknown functions
$\boldsymbol{W}_{j}\big(\boldsymbol{y};\boldsymbol{x},t)\big),\,\,j=ol,\,sp$
(liquid displacements in the variables $\boldsymbol{y}$) we consider the integral identity \eqref{eq2.22} with arbitrary test functions $\displaystyle\,{\xi}=\varepsilon{\eta}(\boldsymbol{x},t)\phi(\frac{\boldsymbol{x}}{\varepsilon})$, where ${\eta}(\boldsymbol{x},t)$ is the same as before and $\phi(\boldsymbol{y})$ is 1-periodic in $\boldsymbol{y}$ function.

Furthermore, using statements~1) and~2) from the conditions of Lemma~\ref{l5.2} we obtain that
\begin{multline}\label{eq5.6}
0=\lim_{\varepsilon\rightarrow{0}}\int_{0}^{{t}_{0}}
\int_{\Omega}\varepsilon{\eta}\phi\chi^{\varepsilon}
(\nabla\cdot\widetilde{\boldsymbol{w}}_{j}^{\varepsilon})dxdt=
\\
\int_{0}^{{t}_{0}}\int_{\Omega}{\eta}\int_{\textbf{Y}_{f}}
\big(\phi\nabla_{y}\cdot\boldsymbol{W}_{j}\big)dydxdt,\,\,j=ol,\,sp.
\end{multline}
Due to the arbitrary choice of the functions ${\eta}$ and $\phi$ the last relation means that hold true the
continuity equation
\begin{equation}\label{eq5.7}
\nabla_{y}\cdot\boldsymbol{W}_{j}(\boldsymbol{y};\boldsymbol{x},t)=0,\,\,\,
(\boldsymbol{y};\boldsymbol{x},t)\in\textbf{Y}_{f}\times\Omega_{t_{0}},\,\,j=ol,\,sp
\end{equation}
together with boundary and normalization conditions
\begin{equation}\label{eq5.8}
\boldsymbol{W}_{j}=0,\,\,(\boldsymbol{y};\boldsymbol{x},t)\in\gamma\times\Omega_{t_{0}},\,\,j=ol,\,sp.
\end{equation}
Let additionally in \eqref{eq2.20} $\boldsymbol{\varphi}=0$ in $\Omega_{T}$, at $t=0$ and $t=t_{0}$,

$\displaystyle\,\frac{\partial\boldsymbol{\varphi}}{\partial{t}}=
\eta(\boldsymbol{x},t)\boldsymbol{\psi}(\frac{\boldsymbol{x}}{\varepsilon})$, where
$\eta\in\mathbb{C}^{\infty}(\overline{\Omega}_{T})$, $\eta(\boldsymbol{x},t)=0$ for $\boldsymbol{x}\in{S}^{0}$, $0<t<T$, $\boldsymbol{\psi}\in{W}^{1}_{2}(Y_{f})$,
\emph{supp}$\boldsymbol{\psi}\subset{Y}_{f}$, $\nabla_{y}\cdot\boldsymbol{\psi}=0$.

Then
\begin{equation*}
\begin{aligned}
& 
\mathbb{D}(x,{\eta}\frac{\partial\boldsymbol{\varphi}}{\partial{t}})=
\mathbb{D}(x,{\eta}\boldsymbol{\psi})=\frac{1}{2}\sum_{i,j=1}^{3}\Big({d}_{ij}\big(x,{\eta}
\boldsymbol{\psi}(\boldsymbol{y})\big)\boldsymbol{e}^{i}\otimes\boldsymbol{e}^{j}+
{d}_{ji}(x,{\eta}\boldsymbol{\psi}(\boldsymbol{y})\big)\boldsymbol{e}^{j}\otimes\boldsymbol{e}^{i}\Big)=
\\
& \frac{1}{2\varepsilon}{\eta}\big(\frac{\partial{\psi}_{j}}{\partial{y}_{i}}(\boldsymbol{y})
\boldsymbol{e}^{i}\otimes\boldsymbol{e}^{j}+\frac{\partial{\psi}_{i}}{\partial{y}_{j}}(\boldsymbol{y})\big)
\boldsymbol{e}^{j}\otimes\boldsymbol{e}^{i}\big)+
\frac{1}{2}\big(\frac{\partial{\eta}}{\partial{x}_{i}}{\psi}_{j}(\boldsymbol{y})
\boldsymbol{e}^{i}\otimes\boldsymbol{e}^{j}+
\frac{\partial{\eta}}{\partial{x}_{j}}{\psi}_{i}(\boldsymbol{y})
\boldsymbol{e}^{j}\otimes\boldsymbol{e}^{i}\big),
\\
& \varepsilon^{2}\mathbb{D}(x,\frac{\partial\boldsymbol{\varphi}}{\partial{t}})=
\varepsilon^{2}\mathbb{D}(x,{\eta}\boldsymbol{\psi})=
\eta\varepsilon\mathbb{D}\big(y,\boldsymbol{\psi}(\frac{\boldsymbol{x}}{\varepsilon})\big)+
\frac{\varepsilon^{2}}{2}(\nabla{\eta}\otimes\boldsymbol{\psi}+
\boldsymbol{\psi}\otimes\nabla{\eta}),
\\
& \nabla\cdot(\eta\boldsymbol{\psi})=
(\nabla{\eta}\cdot\boldsymbol{\psi})=-\frac{1}{\varepsilon}\eta(\nabla_{y}\cdot\boldsymbol{\psi})=0.
\end{aligned}
\end{equation*}
Next we consider functions
\begin{equation}\label{eq5.9}
\boldsymbol{A}_{j}(\boldsymbol{y};\boldsymbol{x},t)=
\nabla\cdot\big(\mu_{j}\mathbb{D}(y,\boldsymbol{W}_{j})\big)\,
\mbox{and}\,\boldsymbol{B}_{j}(\boldsymbol{x},t)=
(t\pi_{j}-{p}^{0})\int_{\textbf{Y}_{f}}\boldsymbol{\psi}dy,\,j=sp,\,ol
\end{equation}
and functional
\begin{equation*}
I^{\varepsilon}_{j}(\eta\boldsymbol{\psi})=\int_{0}^{t_{0}}\int_{\Omega}\chi^{\varepsilon}\eta
\int_{\textbf{Y}_{f}}\Big(\big(\varepsilon^{2}
\widetilde{\boldsymbol{w}}_{f}^{\varepsilon})\cdot\boldsymbol{\psi}\big)+
\big(\mu_{j}\varepsilon^{2}\mathbb{D}(x,\widetilde{\boldsymbol{w}}_{j}^{\varepsilon})+
(t\widetilde{{\pi}}^{\varepsilon}_{j}-{p}^{0})\mathbb{I}\big):\mathbb{D}(x,\boldsymbol{\psi})\Big)dydxdt.
\end{equation*}
In accordance with Lemma~5.2 and the integral identity \eqref{eq2.20} we get
\begin{equation}\label{eq5.10}
\begin{aligned}
& 
0={I}^{0}_{f}(\eta\boldsymbol{\psi})=
\lim_{\varepsilon\rightarrow{0}}I^{\varepsilon}_{f}(\eta\boldsymbol{\psi})=
-\lim_{\varepsilon\rightarrow{0}}\int_{0}^{t_{0}}\int_{\Omega}\chi^{\varepsilon}
\Big((\varepsilon^{2}\widetilde{\boldsymbol{w}}_{f}^{\varepsilon}+\nabla{p}^{0}t\cdot
\frac{\partial\boldsymbol{\varphi}}{\partial{t}})+
\\
& \big(\mu_{j}\varepsilon^{2}\mathbb{D}(x,\widetilde{\boldsymbol{w}}_{j}^{\varepsilon})-
(t\widetilde{{\pi}}^{\varepsilon}-{p}^{0})\mathbb{I}\big):
\mathbb{D}(x,\frac{\partial\boldsymbol{\varphi}}{\partial{t}})\Big)dxdt=
\\
& -\lim_{\varepsilon\rightarrow{0}}\int_{0}^{t_{0}}\int_{\Omega}\eta\chi^{\varepsilon}
\Big(\big(\nabla(\varepsilon\widetilde{\boldsymbol{w}}_{f}^{\varepsilon})\cdot\boldsymbol{\psi}\big)+
\mu_{j}\varepsilon\mathbb{D}(x,\widetilde{\boldsymbol{w}}_{j}^{\varepsilon})\big):
\mathbb{D}\big(y,\boldsymbol{\psi}(\frac{\boldsymbol{x}}{\varepsilon})\big)-
(t\widetilde{{\pi}}^{\varepsilon}-{p}^{0})\big(\nabla\eta\cdot
\boldsymbol{\psi}(\frac{\boldsymbol{x}}{\varepsilon})\big)\Big)dxdt=
\\
& -\int_{0}^{t_{0}}\int_{\Omega}\Big(\eta\big(\int_{Y_{f}}\big(-\nabla_{y}\cdot
\big(\mu_{j}\mathbb{D}(y,\boldsymbol{W}_{j})
\big)dy+(t{\pi}-{p}^{0})(\int_{\textbf{Y}_{f}}\boldsymbol{\psi}dy\cdot\nabla\eta)\Big)dxdt=
\\
& \int_{0}^{t_{0}}\int_{\Omega}\big((B)_{j}<\nabla\eta>-A_{j}\eta\big)dxdt=0.
\end{aligned}
\end{equation}
The last identity in \eqref{eq5.10}
\begin{equation}\label{eq5.11}
\int_{0}^{t_{0}}\int_{\Omega}\big((B)_{j}<\nabla\eta>-A_{j}\eta\big)dxdt=0
\end{equation}
means that functions ${\pi}_{j}\in\mathbb{W}^{1,0}_{2}(\Omega_{T})$ and identity
\eqref{eq5.11} takes the form of the differential equation
\begin{equation}\label{eq5.12}
\nabla_{y}\cdot\big(\mu_{j}\mathbb{D}(y,\boldsymbol{W}_{j})\big)=
-\nabla_{x}(t\pi_{j}-{p}^{0})(\boldsymbol{x},t)=
-\sum_{i=1}^{3}\frac{\partial{}}{\partial{x_{i}}}({\pi}_{j}-{p}^{0}t)(\boldsymbol{x},t)\boldsymbol{e}^{i},
\end{equation}
completed with the continuity equation \eqref{eq5.7}, boundary condition \eqref{eq5.8} and boundary condition \eqref{eq5.5}
\begin{equation}\label{eq5.13}
t\pi_{j}(\boldsymbol{x},t)-{p}^{0}=0,\,\,\,\boldsymbol{x}\in{S}^{1}\cup{S}^{2},\,\,0<t<T,
\end{equation}
which is a consequence of the identity \eqref{eq5.14}.

To solve the periodic boundary value problem \eqref{eq5.7}, \eqref{eq5.8}, \eqref{eq5.12} we use decomposition
\begin{equation}\label{eq5.14}
\boldsymbol{w}_{j}=\int_{\textbf{Y}_{f}}\boldsymbol{W}_{f}\big(\boldsymbol{y};\boldsymbol{x},t\big)dy=
-\frac{1}{\mu_{j}}\sum_{i=1}^{3}\int_{\textbf{Y}_{f}}\boldsymbol{W}_{f}^{(i)}(\boldsymbol{y})dy
\frac{\partial}{\partial{x_{i}}}({t\pi}_{j}-{p}^{0})
(\boldsymbol{x},t),
\end{equation}
where
\begin{equation}\label{eq5.15}
\begin{aligned}
& 
-\nabla\cdot\big(\mu_{j}\mathbb{D}(y,\boldsymbol{W}_{f}^{(i)})\big)=\boldsymbol{e}^{i},\,\,
\nabla\cdot\boldsymbol{W}_{f}^{(i)}=0,\,\,i=1,2,3,\,\,\boldsymbol{y}\in\textbf{Y}_{f},
\\
& \boldsymbol{W}_{f}^{(i)}=0,\,\,i=1,2,3,\,\,\boldsymbol{y}\in\gamma.
\end{aligned}
\end{equation}
The proof of the existence of solutions to the problem \eqref{eq5.15} is standard and follows from energy estimates
\begin{equation*}
\boldsymbol{W}_{f}^{(i)}\in\mathbb{W}^{(1,0}_{2}(\textbf{Y}_{f}),\,\,
\int_{\textbf{Y}_{f}}(|\boldsymbol{W}_{f}^{(i)}|^{2}+|\mathbb{D}(y,\boldsymbol{W}_{f}^{(i)})|^{2})dy
\leqslant\,{M},\,\,i=1,2,3,
\end{equation*}
which are the result of multiplying the equation in \eqref{eq5.15} by $\boldsymbol{W}_{f}^{(i)}$, integration by parts and use of the boundary condition in \eqref{eq5.15} and estimate \eqref{eq2.9} (see Lemma~\ref{l2.10}).

Next, we define the constant matrix $(B)$ as
\begin{equation}\label{eq5.16}
\begin{aligned}
& 
(B)=\sum_{i,j=1}^{3}({b}^{(i,j)}\boldsymbol{e}^{i}\otimes\boldsymbol{e}^{j}+
{b}^{(i,j)}\boldsymbol{e}^{j}\otimes\boldsymbol{e}^{i}),\,\,
\mathbb{B}<\boldsymbol{e}^{i},\boldsymbol{e}^{j}>={b}^{(i,j)},
\\
& {b}^{(i,j)}=
\int_{\textbf{Y}_{f}}(\boldsymbol{W}_{f}^{(i)}(\boldsymbol{y})\cdot\boldsymbol{e}^{j})dy.
\end{aligned}
\end{equation}
The matrix $(B)$ is obviously symmetric and strictly positive definite.
In fact, multiplying the equation in \eqref{eq5.16} by $\boldsymbol{W}_{f}^{(j)}$ gives
\begin{equation*}
\mu_{j}\int_{{Y}_{f}}\mathbb{D}(y,\boldsymbol{W}_{f}^{(i)}):
\mathbb{D}(y,\boldsymbol{W}_{f}^{(j)})dy=
\int_{\textbf{Y}_{f}}(\boldsymbol{W}_{f}^{(j)}\cdot\boldsymbol{e}^{i})dy.
\end{equation*}
Then the equality
\begin{equation*}
\mu_{j}\int_{\textbf{Y}_{f}}\mathbb{D}(y,\boldsymbol{W}_{f}^{(i)}):
\mathbb{D}(y,\boldsymbol{W}_{f}^{(j)})dy=
\int_{\textbf{Y}_{f}}\mu_{j}\mathbb{D}(y,\boldsymbol{W}_{f}^{(j)}):
\mathbb{D}(y,\boldsymbol{W}_{f}^{(i)})dy
\end{equation*}
implies the equality
\begin{equation}\label{eq5.17}
\int_{\textbf{Y}_{f}}(\boldsymbol{W}_{f}^{(i)}\cdot\boldsymbol{e}^{j})dy=
\int_{\textbf{Y}_{f}}(\boldsymbol{W}_{f}^{(j)}\cdot\boldsymbol{e}^{i})dy,
\end{equation}
which means the symmetry of the matrix $(B)_{f}$.

To prove the strict positive definiteness of the matrix $(B)$ we put
\begin{equation*}
\boldsymbol{W}_{f}(\boldsymbol{\xi})=\sum_{i=1}^{3}\boldsymbol{W}_{f}^{(i)}\xi_{i}
\end{equation*}
for any vector $\boldsymbol{\xi}\in\mathbb{R}^{3}$ and any function
$\boldsymbol{\varphi}\in\stackrel{\!\!\circ}{\mathbb{W}}^{1,0}_{2}({Y}_{f})$ and consider
the integral identity
\begin{equation}\label{eq5.18}
\int_{\textbf{Y}_{f}}\mathbb{D}\big(y,\boldsymbol{W}_{f}(\boldsymbol{\xi})\big):
\mathbb{D}(y,\boldsymbol{\varphi})dy=
\int_{\textbf{Y}_{f}}(\boldsymbol{\varphi}\cdot\boldsymbol{\xi})dy.
\end{equation}
Then for $\boldsymbol{\varphi}=\boldsymbol{W}_{f}(\boldsymbol{\xi})$ one has
\begin{equation}\label{eq5.19}
\begin{aligned}
& 
\int_{\textbf{Y}_{f}}|\mathbb{D}(y,\boldsymbol{W}_{f}(\boldsymbol{\xi})|^{2}dy=
\int_{\textbf{Y}_{f}}\big(\boldsymbol{\xi}\cdot\boldsymbol{W}_{f}(\boldsymbol{\xi})\big)dy=
\\
& \sum_{i,j=1}^{3}{\xi}_{i}{\xi}_{j}\int_{\textbf{Y}_{f}}(\boldsymbol{W}_{f}^{(i)}
\cdot\boldsymbol{e}^{j})dy=\sum_{i,j=1}^{3}{b}_{i,j}{\xi}_{i}{\xi}_{j}>0.
\end{aligned}
\end{equation}
It is evident that the equality $\displaystyle\,\sum_{i,j=1}^{3}{b}_{i,j}{\xi}_{i}{\xi}_{j}=0$
implies the equalities $\boldsymbol{W}_{f}=0$ and $\mathbb{D}\big(y,\boldsymbol{W}_{f}(\boldsymbol{\xi})\big)=0$ in $\textbf{Y}_{f}$, which is impossible.

Therefore, we can limit ourselves to the case $|\boldsymbol{\xi}|=1$. This fact immediately leads to the inequality
\begin{equation*}
\sum_{i,j=1}^{3}{b}_{i,j}{\xi}_{i}{\xi}_{j}\geqslant\,\alpha_{0}=\mbox{const}>0.
\end{equation*}
\end{proof}

\section{Conclusion}
The proposed mathematical model at the microscopic level is based on very simple and clear postulates of Newton's classical mechanics of continuous media (section 1.1). These have been verified in practice over several centuries. However, such precise modelling makes any numerical illustration of the physical process impossible. Therefore, we use the homogenization method (section 1.1), proposed by J. Keller and E. Sánchez Palencia, which allows us to move from a microscopic description to a macroscopic description (section 1.1), for which numerical implementations are already possible, visualising all possible variants of fluid flow in real time and at real scale (Lemma 2.12, equations (2.12)–(2.15)).

\begin{funding}
This work was partially supported by the Committee of Science of the Ministry of Science and Higher Education of the Republic of Kazakhstan, grant No. AP32715379 "Mathematical models of filtration of immiscible liquids" (2026-2028).
\end{funding}

\end{document}